\documentclass[12pt, a4paper, reqno, oneside]{article}

\usepackage{amssymb,amsmath}
\DeclareUnicodeCharacter{0335}{}
\usepackage{tikz-cd}
\usepackage{authblk} 
\usepackage{physics} 
\usepackage{tikz}
\usepackage{titlesec} 
\usetikzlibrary{babel}
\usepackage{adjustbox}
\usepackage{url}
\usepackage{xcolor}
\usepackage{yhmath}
\usepackage{eso-pic}
\usepackage[T1]{fontenc}
\usepackage[utf8]{inputenc}
\usepackage[spanish,english]{babel}
\usepackage{csquotes}
\usepackage{pst-node}
\usepackage[english]{babel}
\usepackage{latexsym}
\usepackage{euscript}
\usepackage{mathtools}
\usepackage[all]{xy}
\usepackage{mathrsfs}
\usepackage[percent]{overpic}
\usepackage{setspace}
\usepackage{graphicx}
\usepackage{pdfsync}
\usepackage{amsthm}
\usepackage{thmtools}
\usepackage{yfonts}
\usetikzlibrary{calc}
\usetikzlibrary{decorations.pathmorphing}

\tikzset{curve/.style={settings={#1},to path={(\tikztostart)
    .. controls ($(\tikztostart)!\pv{pos}!(\tikztotarget)!\pv{height}!270:(\tikztotarget)$)
    and ($(\tikztostart)!1-\pv{pos}!(\tikztotarget)!\pv{height}!270:(\tikztotarget)$)
    .. (\tikztotarget)\tikztonodes}},
    settings/.code={\tikzset{quiver/.cd,#1}
        \def\pv##1{\pgfkeysvalueof{/tikz/quiver/##1}}},
    quiver/.cd,pos/.initial=0.35,height/.initial=0}

\tikzset{tail reversed/.code={\pgfsetarrowsstart{tikzcd to}}}
\tikzset{2tail/.code={\pgfsetarrowsstart{Implies[reversed]}}}
\tikzset{2tail reversed/.code={\pgfsetarrowsstart{Implies}}}
\tikzset{no body/.style={/tikz/dash pattern=on 0 off 1mm}}

\usepackage[colorlinks=true,linktocpage=true, pagebackref=false, citecolor=black,linkcolor=black]{hyperref}
\usepackage{orcidlink}
\normalbaselines
\makeatletter

\def\ps@myfancy{\let\@mkboth\markboth
 \def\@evenhead{\vbox{\hsize\textwidth 
 \hbox to \textwidth{\sf\mdseries\thepage 
 \rule[-.6ex]{0mm}{2mm} \hfill\sf\large\leftmark}
 \vskip 1pt \hrule}}
 \def\@oddhead{\vbox{\hsize\textwidth 
 \hbox to \textwidth{{\sf\large\leftmark}
 \rule[-.6ex]{0mm}{2mm} \hfill\sf\mdseries{\thepage}}
 \vskip 1pt \hrule}}}

\def\ps@myfancyplain{
 \def\@evenhead{\vbox{\hsize\textwidth%
 \rule[-.6ex]{0mm}{2mm} \hfill }
 \vskip 1pt \hrule
 \vskip\headsep
 \vskip\textheight
 \vskip1pc
 \hbox to \textwidth{\sf\mdseries\thepage 
 \rule[-6ex]{0mm}{2mm} \hfill }}
 \def\@oddhead{\vbox{\hsize\textwidth 
 \vskip 1pt\hrule
 \vskip\headsep
 \vskip\textheight
 \vskip2pc
 \hbox to \textwidth{\hfill\rule[.4ex]{1pc}{2.5pt}
 \sf\mdseries\thepage}
}}}

\def\ps@myemptyfun{
 \def\@evenhead{\vbox{\hsize\textwidth
 \rule[-.6ex]{0mm}{2mm} \hfill }
 \vskip 1pt 
 \vskip\headsep
 \vskip\textheight
 \vskip1pc
 \hbox to \textwidth{\sf\mdseries\thepage 
 \rule[-0.6ex]{0mm}{2mm} 
 \hfill }}
 \def\@oddhead{\vbox{\hsize\textwidth 
 \vskip 1pt 
 \vskip\headsep
 \vskip\textheight
 \vskip2pc
}}}

\makeatother
\providecommand{\proofname}{Demostraci\'on.}
 {\par\noindent{\it Demostraci\'on. }\nopagebreak\normalsize}%

 {\par\noindent{\it #1. }\nopagebreak\normalsize}%
 {\hfill\linebreak[2]\hspace*{\fill}$\square$\\[-1pt]}
\makeatother

\def\sqbullet{\raise.2ex\hbox{\vrule width 3.5pt height 3.5pt}}

{}

{\em}

\newcounter{substep}
\def\thesubstep{\arabic{substep}}

{\em}

\newcounter{subsubstep}
\def\thesubsubstep{\arabic{subsubstep}}

{\em}

\numberwithin{figure}{section}

{}

\newtheoremstyle{mystyle}
  {}
  {}
  {\itshape}
  {}
  {\sf \bfseries}
  {}
{ }
  {\thmname{#1}\thmnumber{{\textcolor{blue}{\, \hspace{-1mm}#2.}}}\thmnote{ (#3)}}

\theoremstyle{mystyle}
\definecolor{royalblue(web)}{rgb}{0.25, 0.41, 0.88}
\hypersetup{
    colorlinks=true,
    linkcolor=royalblue(web),
    filecolor=magenta,      
    urlcolor=cyan,
    pdftitle={Overleaf Example},
    pdfpagemode=FullScreen,
    }

\titleformat{\section}
  {\normalfont\LARGE\bfseries}{\thesection.}{1em}{}
\titleformat{\subsection}
  {\normalfont\Large\bfseries}{\thesubsection.}{1em}{}
\titleformat{\subsubsection}
  {\normalfont\normalsize\itshape}{\thesubsubsection.}{1em}{}

\newtheorem{Teor}{Theorem}[section]

\newtheorem{Prop}[Teor]{Proposition}
\newtheorem{Coro}[Teor]{Corollary}

\newtheorem{Defi}[Teor]{Definition}

\newtheorem{Lema}[Teor]{Lemma}

 \newcommand{\N}{{\mathbb N}}
\newcommand{\Z}{{\mathbb Z}} 
\newcommand{\R}{{\mathbb R}}
 \newcommand{\C}{{\mathbb C}}

\newcommand{\mail}[1]{\small\href{mailto:#1}{#1}}

\newenvironment{Abstract}
{
\begin{center}
\textbf{Abstract}\\
\vspace{0.25cm}
\begin{minipage}{14.5cm}}
{\footnotesize
\end{minipage}
\end{center}}
\begin{document}


	\begin{center}
		{\huge {\bfseries Fractional weighted Sobolev spaces associated to the Riesz fractional gradient}\par}
		\vspace{1cm}
		
\begin{center} 
	 {\Large Guillermo García-Sáez}{\small\textsuperscript{1}} \\
	\mail{guillermo.garciasaez@uclm.es}\:\orcidlink{0009-0008-4335-119X}
\end{center}
\vspace{5mm}

\textsc{\textsuperscript{1}ETSII, Departamento de Matem\'aticas\\ Universidad de Castilla-La Mancha} \\
		Campus Universitario s/n, 13071 Ciudad Real, Spain. \\ \vspace{5mm}
\end{center}
\begin{Abstract}
In this work, we introduce a new family of functions spaces, the weighted fractional Sobolev spaces $X^{s,p}_{0,w}(\Omega)$, where $w$ is a weight in the Muckenhoupt class $A_p$. This space is a natural extension of the fractional Sobolev spaces $H^{s,p}_0$, obtained by means of the Riesz fractional gradient $\nabla^s$, to the setting of the weighted Lebesgue spaces $L^p_w$. As it happened in the unweighted setting, the spaces $X^{s,p}_{0,w}(\Omega)$ coincide with the weighted version of the Bessel potential space. We obtaien several structural properties for these spaces, as well as continuous and compact embeddings. We conclude with the study of a family of degenerate fractional elliptic partial differential equations.
\end{Abstract}



\tableofcontents

\section{Introduction}
Nowadays, one of the most important topics in the field of calculus of variations and partial differential equations are the nonlocal models arising from nonlocal operators in the form of integro-differential equations. One of the most widely studied nonlocal operators in the past years is the Riesz fractional gradient, which is defined for at least compactly supported smooth functions as $$\nabla ^su=c_{n,s}\int_{\R^n}\frac{u(x)-u(y)}{|x-y|}\frac{x-y}{|x-y|}\,dy.$$ Originally introduced in \cite{ShiehSpector2015}, is the only physically relevant model for a fractional gradient as it was proved in \cite{Silhavy2020} together with its strong convergence to the classical gradient when $s\to 1$ \cite{BellidoCuetoMoraCorral2021}. Further developed in \cite{BellidoCuetoMoraCorral2021, BellidoCuetoGarcia2025, Brue2022, ComiStefani2019, ComiStefani2023, Hidde2022, Scott2022, Silhavy2020}, has encountered many interesting applications to variational models in nonlocal elasticity \cite{Almi2025, BellidoCuetoMoraCorral2020, BellidoCuetoMora2023, Silhavy2024} as well as for more general PDEs problems \cite{Borthagaray2025, camposI,camposII,camposIII,campos2024, Caponi2025, Carrero2025, Gambera2024, Lo2023, Schikorra2015, Schikorra2018, ShiehSpector2018}. Research on the fractional gradient is growing very fast in many different directions, however, as far as the author knows, all the work has been done on the setting of Lebesgue spaces except to the recent work by Campos \cite{campos2024}, in which generalized Orlicz norms are considered to introduce the fractional Orlicz-Sobolev spaces. Those spaces generalize the Orlicz-Sobolev spaces $W^{1,A}_0(\Omega)$ which are used to study nonhomogeneous
elastic materials on the local case. Motivated by this idea and many of the techniques used in the Orlicz framework, we introduce the weighted fractional Sobolev spaces $X^{s,p}_{0,w}(\Omega)$ as the generalization of the spaces $$X^{s,p}_0(\Omega):=\overline{C_c^\infty(\Omega)}^{\norm{\cdot}_{\nabla^s,p}},$$ where $$\norm{u}_{s,p}:=\norm{u}_{L^p}+\norm{\nabla^s u}_{L^p},$$ changing the $L^p$-norm with the weighted $L^p$-norm, i.e., $\norm{u}_{L^p_w}:=\norm{uw^{1/p}}_{L^p}$, where $w$ are nonnegative locally integrable functions with $w>0$ a.e., and such that $1/w$ is locally integrable. In particular we will focus on weights on the $A_p$-class for $1<p<\infty$, which are the weights that verify that the Hardy-Littlewood maximal function is bounded from $L^p_w$ to $L^p_w$. Those kind of weights allow us to use many important results as the $L^p_w$-boundedness of Calderon-Zygmund operators and the H\"ormander-Mihlin multiplier theorem. 

By means of those results we can prove that the spaces $X_{0,w}^{s,p}(\Omega)$ coincide with the weighted version of the Bessel potential spaces $H^{s,p}_{0,w}(\Omega)$, as it happened in the unweighted case \cite[Theorem 1.7]{ShiehSpector2015}. Moreover, our spaces could be obtained by means of complex interpolation as it happened as well for the unweighted case (see \cite{BellidoGarcia2025} and the references therein). Complex interpolation of weighted Sobolev spaces is the main tool on the study of weighted Bessel potential spaces introduced in \cite{Frohlich2004, Schumacher2009I, Schumacher2009II}. Recently, in \cite{Roodenburg2025}, the complex interpolation of weighted Bessel potential spaces is studied for more general weights than the $A_p$-class.

\subsection{Outline of the work}
The structure of the work is as follows: first, we introduce the main definitions and results of complex interpolation, weighted spaces and multiplier theory, and the fractional calculus for $\nabla^s$. Then, in section 3 and 4 we study the weighted fractional Sobolev spaces, establishing the relationship with the weighted Bessel potential spaces and the complex interpolation of weighted Sobolev spaces. Finally, in section 5, we apply those results to study the existence and uniqueness for the problem $$\begin{cases}
    -\operatorname{div}^s\left(w(x)|\nabla^su|^{p-2}\nabla^su\right)&=f,\,x\in \Omega,\\
    \hfill u&=0,\,x\in \Omega^c
\end{cases},$$ which can be seen as the Dirichlet problem for a weighted fractional $p$-Laplacian. For the case $p=2$, we also study the more general problem $$\begin{cases}
    -\operatorname{div}^s\left(A(x)\nabla^su\right)&=f,\,x\in \Omega,\\
    \hfill u&=0,\,x\in \Omega^c
\end{cases},$$ where $A(x)$ is a symmetric $n\times n$ matrix such that $$c_1w(x)|\xi|^2\leq A(x)\xi\cdot \xi\leq c_2w(x)|\xi|^2,$$ for a weight $w\in A_2$. This problem is the fractional counterpart of the prototypical degenerate elliptic partial differential equations. Those kind of equations has been widely studied, in particular in the celebrated work of Fabes, Kenig, and Serapioni \cite{Fabes1982}, where they proved H\"older continuity
and Harnack’s inequality for weak solutions. More recent works on the subject are  \cite{Claros2025, Kebiche2025} in the local case, and in \cite{Behn2024} the nonlocal one.
\section{Preliminaries}
We review the tools that we will require for our development. About the notation, we fix $n\in \mathbb{N}$ the dimension of our ambient space $\R^n$ and we will denote by $\Omega\subset \R^n$ and open bounded subset. The notation for Sobolev $W^{1,p}$ and Lebesgue $L^p$ spaces is the standard one, as is that of smooth
functions of compact support $C_c^\infty$. We will
indicate the domain of the functions, as in $L^p
(\Omega$); the target is indicated only if it is not $\R$.  For two normed spaces $X$, $Y$, we denote the continuous embedding of $X$ into $Y$ as $X\to Y$, i.e., there exists a positive constant $C>0$ such that $\norm{x}_Y\leq C\norm{x}_Y$, for every $x\in X$. For $\alpha \in\mathbb{N}^n$, we give the standard meaning to the partial derivative $\partial^\alpha$ and the size $|\alpha|$.

The ball centered at $x\in \R^n$ and with radius $r>0$ is denoted by $B(x,r)=\{y\in \R^n: |x-y|<r\}$, unless $x=0$, in that case we just write $B_r$. The complementary of a set $E\subset \R^n$ is denoted by $E^c=\R^n\setminus{E}$.

Our convention for the Fourier transform of functions $f\in L^1(\R^n)$ is $$\widehat{f}(\xi)=\int_{\R^n}f(x)e^{-2\pi i x\cdot \xi}\,dx,\,\xi\in \R^n.$$ This definition is extended by continuity and duality to other function and distribution spaces as usually in function spaces theory.
The Schwartz space is denoted by $\mathcal{S}$. The variable
in the Fourier space is generically taken to be $\xi$. We will sometimes use the alternative notation $\mathfrak{F}(f)$ for $\widehat{f}$. More details of this operator could be found in the \cite{Grafakos2008}.
\subsection{Interpolation theory}
In this subsection we briefly briefly introduce Interpolation Theory, including the concepts and results we will use in the following. For a complete development of the main ideas, we refer to \cite{BerghLofstrom1976,GarciaSaez2024,Lunardi2018,Triebel1995}.\\
Let $(E_0,E_1)$ a couple of Banach spaces. We say that the couple is \textit{compatible} if there exists a Hausdorff topological vector space $\mathcal{E}$ such that $E_0,E_1\xhookrightarrow{}\mathcal{E}$. We say that a Banach space $E$ is intermediate with respect to the couple if $E_0\cap E_1\xhookrightarrow{}E\xhookrightarrow{}E_0+E_1$. Let $(F_0,F_1)$ another compatible couple and $T:E_0+E_1\to F_0+F_1$ bounded a linear. We say that $T$ is \textit{admissible} if $T:E_j\to F_j$, $j=0,1,$ continuously. Given $E$ an intermediate space with respect to the couple $(E_0,E_1)$, and $F$ an intermediate space with respect to $(F_0,F_1)$, we say that $E,F$ are \textit{interpolation spaces} with respect to the couples $(E_0,E_1)$ and $(F_0,F_1)$, respectively, if for every admissible operator $T:E_0+E_1\to F_0+F_1$, $T:E\to F$ continuously. The methods to construct such interpolation spaces for given couples are called \textit{interpolation functors.} We say that an interpolation functor $\mathcal{F}$ is of \textit{exponent} $\theta\in (0,1)$ if $$\norm{T}_{\mathcal{F}\left((E_0,E_1)\right)\to \mathcal{F}\left((F_0,F_1)\right)}\leq C\norm{T}_{E_0\to F_0}^{1-\theta}\norm{T}_{E_1\to F_1}^\theta,$$ for some positive constant $C$. If we can choose $C=1$, we say that the functor is \textit{exact} of exponent $\theta$.\\

For our purposes, we will focus on the complex interpolation method. Given a compatible couple of Banach spaces $(E_0,E_1)$, we define the space $\mathfrak{F}(E_0,E_1)$, as the space of functions $f:S\to E_0+E_1$, where $S:=\{z\in \C: 0\leq \operatorname{Re}z\leq 1\}$, such that $f$ is holomorphic on the interior of $S$, continuous and bounded on $S$, and the functions $t\mapsto f(j+it)$, $j=0,1$, are continuous from $\R\to E_j$, and such that $\norm{f(j+it)}_{E_j}\to 0$ as $|t|\to \infty$. The space $\mathfrak{F}(E_0,E_1)$ is a vector space which is complete endowed with the norm $$\norm{f}_{\mathfrak{F}(\overline{E})}:=\operatorname{max}\{\operatorname{sup}_{t\in\R}\norm{f(it)}_{E_0},\operatorname{sup}_{t\in\R}\norm{f(1+it)}_{E_1}\},\,f\in \mathfrak{F}(\overline{E}).$$ From this space, we construct the \textit{complex method} as the functor $\mathcal{C}_\theta$, $\theta\in [0,1]$ which associates the space $[E_0,E_1]_\theta$ to the compatible couple of Banach spaces $(E_0,E_1)$. The space $[E_0,E_1]_\theta$ is defined as the space of $x\in E_0+E_1$ such that there exists $f\in \mathfrak{F}(E_0,E_1)$ with $f(\theta)=x$. The space is a Banach space endowed with the norm $$\norm{x}_{\theta}:=\operatorname{inf}\{\norm{f}_{\mathfrak{F}(\overline{E})}: f\in \mathfrak{F}(\overline{E}), f(\theta)=x\}.$$ 
The spaces $[E_0,E_1]_\theta$ are exact interpolation spaces of exponent $\theta$ (see \cite[Theorem~4.1.2]{BerghLofstrom1976} or \cite[Theorem~IV.1.5]{GarciaSaez2024}).\\
Now we enumerate the most important properties of such spaces.
\begin{Teor}[Properties of Complex interpolation spaces]\label{ComplexProps}
        Let $(E_0,E_1)$ a compatible couple of Banach spaces, and $\theta\in [0,1]$. Then, we have
        \begin{enumerate}
            \item[1.-] $[E_0,E_1]_{\theta}=[E_1,E_0]_{1-\theta}$.
            \item[2.-]If $E_0\xhookrightarrow{}E_1$ and $\theta_0<\theta_1$, $[E_0,E_1]_{\theta_0}\xhookrightarrow{}[E_0,E_1]_{\theta_1}$.
            \item[3.-] If $E_0=E_1$ and $0<\theta<1$, $[E_0,E_1]_\theta=E_0$.
            \item[4.-] $E_0\cap E_1$ is dense in $[E_0,E_1]_\theta$.
            \item[5.-] There exists a positive constant $C>0$ such that for every $u\in E_0\cap E_1$, $$\norm{u}_{[E_0,E_1]_\theta}\leq C\norm{u}_{E_0}^{1-\theta}\norm{u}_{E_1}^\theta.$$
            \item[5.-] $[E_0,E_1]_j$ is a closed subspace of $E_j$ with coincidence of the norm in $[E_0,E_1]_j$, $j=0,1$.
            \item[6.-] If $E_0\cap E_1$ is dense in both $E_0$ and $E_1$, and at least one of $E_0$ or $E_1$ is reflexive, then $$\left([E_0,E_1]_\theta\right)^*=[E_0^*,E_1^*]_\theta,\,\theta\in (0,1).$$
            \item[7.-]  If at least one of $E_0$ or $E_1$ is reflexive, then the space $[E_0,E_1]_\theta,\,\theta\in (0,1),$ is reflexive
            \item[8.-](Complex reiteration theorem) For every $\theta_0,\theta_1\in (0,1)$ and $\alpha\in (0,1)$, $$\left[[E_0,E_1]_{\theta_0},[E_0,E_1]_{\theta_1}\right]_\alpha=[E_0,E_1]_{\theta(\alpha)},\,\theta(\alpha)=(1-\alpha)\theta_0+\alpha\theta_1.$$ Moreover, for $1\leq q\leq \infty$, $$\left([E_0,E_1]_{\theta_0},[E_0,E_1]_{\theta_1}\right)_{\alpha,q}=(E_0,E_1)_{\theta(\alpha),q}.$$
         \end{enumerate}
\end{Teor}
Detailed proofs of these facts can be found in \cite[Theorem~4.2.1, Theorem~4.2.2, Theorem~4.5.1]{BerghLofstrom1976} and \cite[Proposition~IV.1.8, Theorem~IV.5.4, Theorem~IV.5.6]{GarciaSaez2024}.
\subsection{Weighted spaces}
By a cube $Q$ in $\R^n$ we mean a subset of $\R^n$ on the form $I_1\times \cdots\times I_n$, where each $I_j$, $j=1,\ldots,n$ is a bounded interval of $\R^n$, each one of the same length. The cubes have always sides parallel to the axes.

A \textit{weigth} $w$ is a nonnegative locally integrable function on $\R^n$ taking values in $(0,\infty)$ almost everywhere, i.e., $w$ is allowed to be zero or infinity only on a set of zero Lebesgue measure. Hence, if $w$ is a weight and $1/w$ is locally integrable, then $1/w$ is also a weight. Since we will be dealing with weights on functions spaces derived from $L^p$-norms, we are interested in the following class of weights.
\begin{Defi}\em
    Let $1<p<\infty$. A function $0\leq w\in L_{\text{loc}}^1(\R^n)$ is called an $A_p$-weight if $$[w]_p:=\operatorname{sup}_Q\left(\frac{1}{|Q|}\int_Qw(x)\,dx\right)\left(\frac{1}{|Q|}\int_{Q}w(x)^{-1/(p-1)}\,dx\right)^{p-1}<\infty,$$ where the supremum is taken over all cubes of $\R^n$ and $|Q|$ is the usual Lebesgue measure of the cube. We will use the abreviation $w(U)$ for $\int_U w$, and set $w^*:=w^{-1/(p-1)}$. \end{Defi}
    Note that it is posible to define the $A_p$-class of weights for $p=1$ and $p=\infty$. However, those cases are significantly different and they are not of our interest for this work.
The following Lemma summarizes some important properties of the $A_p$-class for $1<p<\infty$. See \cite[Ch. 4]{garcia1985weighted} for detailed proofs.
\begin{Lema}
    Let $1<p,q<\infty$. Then, \begin{itemize}
        \item $w\in A_p\iff w^*\in A_{p'}$, with $1/p+1/p'=1$.
        \item $A_p\subset A_q$ if $p\leq q$.
        \item  For $w\in A_p$, there exists $r<p$ such that $w\in A_r$.
    \end{itemize}
\end{Lema}   
 Let $1<p<\infty$, $w$ a weight and $\Omega$ a Lipschitz bounded domain on $\R^n$ or $\Omega=\R^n$. We define the \textit{weighted $L^p$ spaces} as $$L^p_w(\Omega):=\Bigg\{u\in L_{\operatorname{loc}}^1\left(\overline{\Omega}\right): \norm{u}_{L^p_w(\Omega)}:=\left(\int_\Omega |f|^pw\,dx\right)^{1/p}<\infty\Bigg\}.$$ The condition $w\in A_p$, sometimes referred as the \textit{Muckenhoupt condition}, comes from the natural question of whether the Hardy-Littlewood maximal function of $f\in L^p_w(\Omega)$ is bounded on $L^p_w(\Omega)$, i.e., determining the class of functions $w$ such that the operator $$Mf(x):=\operatorname{sup}_{r>0}\frac{1}{r^n}\int_{B(x,r)}|f(y)|\,dy,$$ is bounded from $L^p_w\to L^p_w$.
 The answer to this question was given by Muckenhoupt in its work \cite{Muckenhoupt}, where it was established that the boundedness of $M$ on $L^p_w$ is equivalent to having $w\in A_p$ for $1<p<\infty$.\\

From now, we are only going to deal with $A_p$-weights.
Many properties known for the usual Lebesgue spaces also holds for the weighted case, for example the dual space characterization for $1<p<\infty$, $$\left(L^p_w(\Omega)\right)^*=L^q_{w^*}(\Omega),$$ where $1/p+1/q=1$. The pairing $\langle f,u\rangle$ for $f\in L^q_{w^*}(\Omega)$ and $u\in L^p_w(\Omega)$ is given by $\int_\Omega fu$.\\

A prototypical example of weight on the $A_p$ class are radially symmetric weights of the form $w(x)=|x-x_0|^\alpha$ for $-n<\alpha<n(p-1)$. Moreover, if we consider the distance function $w(x)=d(x,M)^\alpha$ where $M$ is a compact Lipschitzian manifold of dimension $k<n$, $w\in A_p\iff -(n-k)<\alpha<(n-k)(p-1)$. 

In analogy to the unweighted case, we can also define a weighted version of the classical Sobolev spaces just changing the usual $p$-norm with the $L^p_w$-norm. Let $k$ a natural number, $1<p<\infty$, $w\in A_p$ and $\Omega$ a Lipschitz domain of $\R^n$ or $\Omega=\R^n$. Then, we define the \textit{weighted Sobolev spaces} $W^{k,p}_w(\Omega)$ as the space of $L^p_w(\Omega)$ functions $u$ such that $D^\alpha u\in L^p_w(\Omega)$ for every multi-index $\alpha$ with $|\alpha|\leq k$. The space is endowed with the norm  $$\norm{u}_{W^{k,p}_w(\Omega)}:=\left(\sum_{|\alpha|\leq k}\norm{D^\alpha u}_{L^p_w(\Omega)}\right)^{1/p}.$$
\subsection{Calder\'on-Zymgund operators and multiplier theory}
An important class of operators given by singular integral are the Calder\'on-Zymgund operators. They are very related to the theory of weighted Lebesgue spaces and the will play a key role in our work. We recall the main definitions and some crucial results.\\

We say that a linear operator $T$ is a \textit{Calder\'on-Zygmund operator} if there exists a kernel $K:\R^n\times \R^n\setminus\{(x,x):x\in \R^n\}\to \R$, such that for all $\varphi\in C_c^\infty(\R^n)$ and $x\not \in \operatorname{supp}f$, $$T\varphi(x)=\int_{\R^n}K(x,y)\varphi(y)\,dy,$$ with the kernel satisfying that $$|K(x,y)|\leq C|x-y|^{-n},$$ and that for some $\varepsilon>0$, $$|K(x+h,y)-K(x,y)|+|K(x,y)-K(x,y+h)|\leq C\frac{|h|^\varepsilon}{|x-y|^{n+\varepsilon}},$$ for every $h\in \R^n$ such that $2|h|\leq |x-y|$.
This type of operators satisfy that for every $1<p<\infty$ and $w\in A_p$, $T:L_w^p(\R^n)\to L^p_w(\R^n)$ continuously (see \cite[Theorem 7.11]{duandikoetxea2001}). Indeed, we have the following result concerning the weighted $L^p$-boundedness of many important operators which also characterices the class $A_p$ (see \cite[Theorem 1.1]{Cruz-uribe2011}).
\begin{Teor}
    Let $1<p<\infty$, $w$ a weight and let $T$ denote the Hardy-Littlewood maximal operator, the Hilbert transform or the Riesz transform. Then, $w\in A_p$ if and only if $$\norm{Tf}_{L^p_w(\R^n)}\leq C\norm{f}_{L^p_w(\R^n)},\,f\in L_w^p(\R^n),$$ for some positive constant $C$ no depending on $f$. More generally, this result holds for any operator given as a convolution with a sufficiently smooth kernel.
\end{Teor}
The sufficiency of the $A_p$ weights to the boundedness of Calder\'on-Zymgund operators, which is the result which we will require, is due to Coiffman and Fefferman 
\cite{coifmanfefferman1974}.\\
 A standard procedure to prove that an operator $T$ is a Calder\'on-Zygmund
singular integral operator is by means of the so called H\"ormander-Mihlin multiplier theorem.
\begin{Teor}\label{Milhinnormal}{\normalfont \bfseries {[H\"ormander-Mihlin Multiplier Theorem]}}
    Let $m\in C^n\left(\R^n\setminus\{0\}\right)$ such that $$\bigg|\partial^\alpha_\xi m(\xi)\bigg|\leq C|\xi|^{-|\alpha|},\,\xi\in \R^n,\,\alpha\in\N:|\alpha|=0,1,\ldots,\left[\frac{n}{2}\right],$$ for some positive constant $C>0$. Then, for any $1<p<\infty$, we have that the convolution operator $$T\varphi:=\check{m}*\varphi,\,\varphi\in \mathcal{S}(\R^n),$$ extends to a bounded operator from $L^p(\Omega)\to L^p(\Omega)$, where $\Omega\subseteq \R^n$ is an open set. Here, $\check{m}$ denotes the inverse Fourier transform of $m$.
\end{Teor}
This multiplier theorem also holds for weighted Lebesgue spaces, as it was proven by Garc\'ia-Cueva and Rubio de Francia in \cite[Theorem 3.9, pp.418]{garcia1985weighted}:
\begin{Teor}\label{Milhin}{\normalfont \bfseries {[H\"ormander-Mihlin Multiplier Theorem with Weights]}}
    Let $m\in C^n\left(\R^n\setminus\{0\}\right)$ such that $$\bigg|\partial^\alpha_\xi m(\xi)\bigg|\leq C|\xi|^{-|\alpha|},\,\xi\in \R^n,\,\alpha\in \N:|\alpha|=0,1,\ldots,\left[\frac{n}{2}\right]+1,$$ for some positive constant $C>0$. Then, for any $1<p<\infty$ and $w\in A_p$, we have that the convolution operator $$T\varphi:=\check{m}*\varphi,\,\varphi\in \mathcal{S}(\R^n),$$ extends to a bounded operator from $L^p_w(\Omega)\to L^p_w(\Omega)$, where $\Omega\subseteq \R^n$ is an open set.
\end{Teor}
\subsection{Nonlocal calculus}
Since we will be dealing with weighted fractional Sobolev spaces, we recall the basic definitions of the fractional operators introduced by Shieh and Spector in \cite{ShiehSpector2015,ShiehSpector2018}. Let $u\in C_c^\infty(\R^n)$, $v\in C_c^\infty(\R^n;\R^n)$ and $s\in (0,1)$. We define the \textit{Riesz fractional gradient} $\nabla^s$ as $$\nabla^su(x)=c_{n,s}\int_{\R^n}\frac{u(x)-u(y)}{|x-y|^{n+s}}\frac{x-y}{|x-y}\,dy,\,x\in \R^n,$$ and the \textit{fractional divergence} as $$\operatorname{div}^sv(x)=c_{n,s}\int_{\R^n}\frac{v(x)-v(y)}{|x-y|^{n+s}}\cdot\frac{x-y}{|x-y}\,dy,\,x\in \R^n,$$ where $$c_{n,s}=\frac{\Gamma\left(\frac{n+s+1}{2}\right)}{\pi^{n/2}2^{-s}\Gamma\left(\frac{1-s}{2}\right)},$$ is a normalizing constant. One of the main properties of these objects is that they can be seen as the classical gradient and divergence of the Riesz potential $$I_su:=I_{s}*u(x),$$ where $$I_s(x):=\frac{1}{\gamma_{n,s}}\frac{1}{|x|^{n-s}},\,0<s<n.$$  of the functions $u$ and $v$. Here, the constant $\gamma_{n,s}:=\frac{\pi^{n/2}2^s\Gamma(s/2)}{\Gamma\left(\frac{n-s}{2}\right)}$, satisfies that $$\gamma_{n,1-s}c_{n,s}=(n+1-s).$$ In particular we have the following result: \begin{Prop}Let $s\in (0,1)$, $u\in C_c^\infty(\R^n)$ and $v\in C_c^\infty(\R^n;\R^n).$ Then, \begin{align*}
    \nabla^su&=D(I_{1-s}u)=I_{1-s}(Du),\\
    \operatorname{div}^sv&=\operatorname{div}(I_{1-s}v)=I_{1-s}(\operatorname{div}v).
\end{align*}
\end{Prop}
\noindent{}\textbf{Proof:} Let $u\in C_c^\infty(\R^n)$ and we fix $\varepsilon>0$. By the first Green's identity, \begin{align*}
    &\int_{B(x,\varepsilon)^c}\frac{D_yu(y)}{|x-y|^{n+s-1}}\,dy=\int_{B(0,\varepsilon)^c}\frac{D_zu(x+z)}{|z|^{n+s-1}}\,dz\\
    &=\int_{\partial B(0,\varepsilon)^c}\frac{u(x+z)}{|z|^{n+s-1}}\cdot \frac{z}{|z|}\,dS_z+\int_{B(0,\varepsilon)^c}u(x+z)D_z\left(\frac{1}{|z|^{n+s-1}}\right)\,dz\\
    &=\int_{\partial B(0,\varepsilon)^c}\frac{u(x+z)}{|z|^{n+s-1}}\cdot \frac{z}{|z|}\,dS_z+(n+s-1)\int_{B(0,\varepsilon)^c}\frac{u(x+z)}{|z|^{n+s+1}}z\,dz.
\end{align*}
By the odd radial symmetry, \begin{align*}\int_{B(0,\varepsilon)^c}\frac{u(x+z)}{|z|^{n+s}}\frac{z}{|z|}\,dz&=\int_{B(x,\varepsilon)^c}\frac{-u(x)+u(x+z)}{|z|^{n+s}}\frac{z}{|z|}\,dz\\
&=\int_{B(x,\varepsilon)^c}\frac{u(x)-u(y)}{|x-y|^{n+s}}\frac{x-y}{|x-y|}\,dy.\end{align*} On the other hand, \begin{align*}
    \Bigg|\int_{\partial B(0,\varepsilon)^c}\frac{u(x+z)}{|z|^{n+s-1}}\cdot \frac{z}{|z|}\,dS_z\Bigg|
    &\leq \norm{Du}_\infty \int_{\partial B(0,\varepsilon)^c}\frac{1}{|z|^{s-1}}\,dS_z\\
    &=\norm{Du}_\infty \frac{2\pi^{n/2}}{\Gamma(n/2)}\varepsilon^{1-s}\to 0,
\end{align*} if $\varepsilon\to 0^+.$ Hence, \begin{align*}
   &I_{1-s}(Du)=\frac{1}{\gamma_{n,1-s}}\int_{\R^n}\frac{D_yu(y)}{|x-y|^{n+s-1}}\,dy=\frac{1}{\gamma_{n,1-s}}\lim_{\varepsilon\to 0^+}\int_{B(x,\varepsilon)^c}\frac{D_yu(y)}{|x-y|^{n+s-1}}\,dy\\
    &=\frac{1}{\gamma_{n,1-s}}\lim_{\varepsilon\to 0^+}\left(\int_{\partial B(0,\varepsilon)^c}\frac{u(x+z)}{|z|^{n+s-1}}\cdot \frac{z}{|z|}\,dS_z+(n+s-1)\int_{B(x,\varepsilon)^c}\frac{u(x)-u(y)}{|x-y|^{n+s}}\frac{x-y}{|x-y|}\,dy\right)\\
    &=\frac{(n+s-1)}{\gamma_{n,1-s}}\lim_{\varepsilon\to 0^+}\int_{B(x,\varepsilon)^c}\frac{u(x)-u(y)}{|x-y|^{n+s}}\frac{x-y}{|x-y|}\,dy\\
    &=\frac{(n+s-1)}{c_{n,s}\gamma_{n,1-s}}c_{n,s}\operatorname{pv}_x\,\int_{\R^n}\frac{u(x)-u(y)}{|x-y|^{n+s}}\frac{x-y}{|x-y|}\,dy,\\
    &=D^su.
\end{align*} The proof for $\operatorname{div}^s$ is completely analogous. \qed 

The following proposition collect some of the most important facts about those operators.
\begin{Prop}\label{DsProps}
    Let $s\in (0,1)$, $u\in C_c^\infty(\R^n)$ and $v\in C_c^\infty(\R^n;\R^n)$.
    \begin{itemize}
    \item The Fourier transform of $\nabla^su$ is $$\widehat{\nabla^su}(\xi)=\frac{2\pi i\xi}{|2\pi\xi|^{1-s}}\widehat{u}(\xi),$$ and the Fourier transform of $\operatorname{div}^su$ is $$\widehat{\operatorname{div}^su}(\xi)=\frac{2\pi i\xi}{|2\pi\xi|^{1-s}}\cdot \widehat{v}(\xi).$$
        \item The following fractional integration by parts formula holds: $$\int_{\R^n}\nabla^su(x)\cdot v(x)\,dx=-\int_{\R^n}u(x)\operatorname{div}^sv(x)\,dx.$$
        \item We have the following fractional theorem of calculus: $$u(x)=c_{n,-s}\int_{\R^n}\nabla^su(y)\cdot \frac{x-y}{|x-y|^{n-s+1}}\,dy=I_s(\mathcal{R}\cdot \nabla^su),\,x\in\R^n,$$ where $\mathcal{R}$ is the Riesz transform, defined as a Fourier multiplier with symbol $-i\frac{\xi}{|\xi|}$.
        \end{itemize}
\end{Prop}
For detailed proofs of these facts and a complete exposition on the fractional setting we refer to \cite{ShiehSpector2015}.\\

\section{Weighted Bessel potential spaces}
For simplicity, we will use the notation $\left\langle\xi\right\rangle:=\left(1+4\pi^2|\xi|^2\right)^{1/2}, \xi\in \R^n$, sometimes called \textit{Japanese bracket}. \begin{Defi}\em
    Let $s\in \C$, $f\in \mathcal{S}'(\R^n)$ (where $ \mathcal{S}'(\R^n)$ is the space of tempered distributions) and $\xi\in \R^n$. We define the \textit{Bessel potential} $\Lambda_s$ of order $s$ of $f$ as $$\Lambda_{s}f=\mathfrak{F}^{-1}\left(\left\langle\xi\right\rangle^{-s}\widehat{f}(\xi)\right).$$
\end{Defi}
The name Bessel potential comes from the fact that for every $s>0$, $$\Lambda_sf=G_s*f,\,f\in \mathcal{S}'(\R^n),$$ where $G_s$ is the \textit{Bessel kernel} defined as $$G_s(x)=\frac{1}{(4\pi)^{n/2}\Gamma(n/2)}\int_0^\infty e^{-t/(4\pi)}e^{-|x|^2\pi/t}t^{(s-n)/2}\,\frac{dt}{t},\,x\in \R^n.$$ 
\begin{Defi}\em[Weighted Bessel potential space] Let $s\in\R$, $1<p<\infty$ and $w\in A_p$. We define the \textit{weighted Bessel potential space} $H^{s,p}_w(\R^n)$ as $$H^{s,p}_w(\R^n):=\{u\in \mathcal{S}'(\R^n): \Lambda_{-s}u\in L^p_w(\R^n)\},$$ endowed with the norm $\norm{u}_{H^{s,p}_w(\R^n)}:=\norm{\Lambda_{-s}u}_{L^p_w(\R^n)}$.
\end{Defi}
We first want to note that those space are well defined since $H^{s,p}_w(\R^n)\xhookrightarrow{}L^p_w(\R^n)$, for $w\in A_p$. This follows from the fact that $\Lambda_s$ maps $L^p_w(\R^n)$ to $L^p_w(\R^n)$ continuously. In fact, since $G_s\in L^1(\R^n)$ and is radial function, by \cite[Theorem 2, pp. 62]{Stein1970}, $G_s*f(x)\leq Mf(x)$, for locally integrable functions $f$, and hence  $$\norm{G_s*f}_{L^p_w(\R^n)}\leq \norm{Mf}_{L^p_w(\R^n)}\leq C\norm{f}_{L^p_w(\R^n)},$$ so for $u\in H^{s,p}_w(\R^n)$, $$\|u\|_{L^p_w(\R^n)}= \|\Lambda_{s}\Lambda_{-s} u\|_{L^p_w(\R^n)}\leq C\|\Lambda_{-s}u\|_{L^p_w(\R^n)}=C\|u\|_{H_w^{s,p}(\R^n)}.$$
\textbf{Theorem \ref{Milhin}}, allows us to establish the interpolation nature of those spaces, as it happened in the unweighted case. In particular, from \textbf{Theorem \ref{Milhin}} it is straightforward to obtain the boundedness of the purely imaginary powers of the Bessel potential on the space $L^p_w(\R^n)$, analogous to \cite[Lemma 3.8]{BellidoGarcia2025}, and hence replicating the proofs of \cite[Theorem 3.7, Proposition 3.9]{BellidoGarcia2025}, we obtain the following result:
\begin{Teor}
    Let $s\in \R$ such that $s=\theta k$, with $\theta\in (0,1)$ and $k\in \mathbb{Z}$, $1<p<\infty$ and $w\in A_p$. Then, $$H^{s,p}_w(\R^n)=[L^p_w(\R^n),W^{k,p}_w(\R^n)]_\theta,$$ with equivalence of the norms.
\end{Teor}
From the main properties of the complex interpolation functor of Theorem \ref{ComplexProps}, the following results about $H^{s,p}_w(\R^n)$ are straightforward.
\begin{Prop}
    Let $s\in \R$, $1<p<\infty$ and $w\in A_p$. Then, \begin{enumerate}
        \item $H_w^{s,p}(\R^n)$ is a complete, reflexive and separable space.
        \item $\left(H^{s,p}_w(\R^n)\right)^*=H^{-s,p'}_{w^*}(\R^n)$, where $1/p+1/p'=1$.
        \item $H^{s,p}_w(\R^n)\xhookrightarrow{}H^{t,p}_w(\R^n)$ for every $0<t<s$.
        \item If $s=k\in \mathbb{N}$, $H^{k,p}_w(\R^n)=W^{k,p}_w(\R^n)$.
        \item $W^{k,p}_w(\R^n)$ is dense in $H^{\theta k,p}_w(\R^n)$ for every $\theta\in (0,1)$ and $k\in \mathbb{N}$. The space $C_c^\infty(\R^n)$ is dense in $H^{s,p}_w(\R^n)$ if $s>0$.
    \end{enumerate}
\end{Prop}

Let $s\in (0,1)$ and $1<p<\infty$. We define the norm $$\norm{u}_{\nabla^s,p}:=\norm{u}_p+\norm{\nabla^s u}_p,\,u\in C_c^\infty(\R^n).$$ We know that $H^{s,p}(\R^n)=\overline{C_c^\infty(\R^n)}^{\norm{\cdot}_{\nabla^s,p}}$ (see \cite[Theorem 1.7]{ShiehSpector2015}). Motivated by this fact, we consider the norm $$\norm{u}_{\nabla^s,p,w}:=\norm{u}_{L_w^p(\R^n)}+\norm{\nabla^su}_{L^p_w(\R^n;\R^n)},\,u\in C_c^\infty(\R^n),$$ where $w\in A_p$. From this norm, we introduce the \textit{fractional weighted Sobolev spaces} as $$X^{s,p}_w(\R^n)=\overline{C_c^\infty(\R^n)}^{\norm{\cdot}_{\nabla^s,p,w}}.$$
To be consequent with our definition of the space, we have to extend the definition of $\nabla^su$ for functions $u\in X_w^{s,p}(\R^n)$ by continuity. Since by definition, for every $u\in X_w^{s,p}(\R^n)$ there exists a sequence $\{u_m\}\subset C_c^\infty(\R^n)$ such that $u_m\to u$ in $L^p_w(\R^n)$ and $\{\nabla^su_m\}$ is Cauchy in $L^p_w(\R^n;\R^n)$, we define $$\nabla^su:=L^p_w(\R^n;\R^n)-\lim_{m\to \infty} \nabla^su_m.$$ From this characterization is it clear that the duality between the operators $\nabla^s$ and $\operatorname{div}^s$, i.e., the integration by parts formula, also holds for functions in $X^{s,p}_w(\R^n)$. This fact implies that the definition of $\nabla^su$ is independent of the choice of the Cauchy sequence $\{u_n\}$, and hence the space $X^{s,p}_w(\R^n)$ is well defined. Moreover, $I_{1-s}$ can be seen as a translation operator from $X^{s,p}_w(\R^n)\to W^{1,p}_w(\R^n)$ as in the unweighted case.\\
As a direct application of the weighted H\"ormander-Mihlin theorem, we get that our spaces coincide with the weighted Bessel potential spaces, as it happened in the unweighted case.
For the proof, we will work with an alternative but equivalent definition of the Riesz fractional gradient in a distributional sense. Let $u\in L^p(\R^n)$ such that $I_{1-s}u$ is well defined. Then, we define the \textit{distributional Riesz fractional gradient} $D^su=(D^su)_j$ as $$(D^su)_j:=\frac{\partial^su}{\partial x_j^s}:=\frac{\partial}{\partial x_j}I_{1-s}u,\,j=1,\ldots,n,$$ in the sense that $$\left\langle \frac{\partial^su}{\partial x_j^s},\varphi\right\rangle:=-\int_{\R^n}I_{1-s}u\frac{\partial \varphi}{\partial x_j}\,dx=-\left\langle I_{1-s}u,\frac{\partial\varphi}{\partial x_j}\right\rangle,$$ for every $\varphi\in C_c^\infty(\R^n)$. Since by \cite[Theorem 1.2]{ShiehSpector2015}, $D^su=I_{1-s}(Du)=D(I_{1-s}u)$, $D^s$ and $\nabla^s$ coincide at least for smooth compactly supported functions.
\begin{Teor}\label{Equivalence}
    Let $s\in (0,1)$, $1<p<\infty$ and $w\in A_p$. Then, $$X^{s,p}_w(\R^n)=H_w^{s,p}(\R^n),$$ with equivalence of the norms.
\end{Teor}
\noindent{}\textbf{Proof:} We first prove the inclusion $X_w^{s,p}(\R^n)\xhookleftarrow{}H_w^{s,p}(\R^n)$. Let $u\in H_w^{s,p}(\R^n)$, so there exists a function $v\in L_w^p(\R^n)$ such that $u=\Lambda^s v$. Since $\Lambda_s$ is a bounded linear operator from $L_w^p(\R^n)$ to $L_w^p(\R^n)$ for $s>0$ (\cite[pp.~75]{turesson2000}), $u\in L_w^p(\R^n)$, so we only need to prove that $\nabla^su\in L_w^p(\R^n;\R^n)$. By the density of $C_c^\infty(\R^n)$ in $L_w^p(\R^n)$, it is enough to consider that $v\in C_c^\infty(\R^n)$. Now, observing that for $\psi\in C_c^\infty(\R^n),$ \begin{align*}\left\langle \frac{\partial^su}{\partial x_j^s},\psi\right\rangle&=-\left\langle I_{1-s}u,\frac{\partial \psi}{\partial x_j }\right\rangle=-\left\langle I_{1-s}\Lambda_s v,\frac{\partial \psi}{\partial x_j }\right\rangle=-\left\langle \Lambda_s\left(I_{1-s}u\right),\frac{\partial \psi}{\partial x_j }\right\rangle\\
&=\left\langle \Lambda_s\left( \frac{\partial^s v}{\partial x_j^s}\right),\psi\right\rangle,\end{align*} we get that $$\frac{\partial^su}{\partial x_j^s}=\Lambda_s\left(\frac{\partial^sv}{\partial x_j^s}\right),$$ and hence \begin{align*}\mathfrak{F}\Bigg\{\frac{\partial^su}{\partial x_j^s}\Bigg\}(\xi)&=\mathfrak{F}\Bigg\{\Lambda_s\left( \frac{\partial^s v}{\partial x_j^s}\right)\Bigg\}(\xi)=(1+4\pi|\xi|^2)^{-s/2}\left(\frac{2\pi i\xi_j}{|2\pi \xi|^{1-s}}\right)\mathfrak{F}\{v\}(\xi)\\
&=\frac{i \xi_j}{|\xi|}\frac{(2\pi|\xi|)^s}{(1+4\pi|\xi|^2)^{s/2}}\mathfrak{F}\{v\}(\xi).\end{align*} 
We define the operator $$\widehat{T_s\varphi}:=\frac{(2\pi|\xi|)^s}{\left\langle\xi\right\rangle^{s}}\widehat{\varphi},\,\varphi\in \mathcal{S}(\R^n).$$
By Theorem \ref{Milhin}, $T_s$ extends to a bounded operator from $L^p_w\to L^p_w$, and by the $L_w^p$-boundedness of the Riesz transform, $$\norm{\frac{\partial^s u}{\partial x_j^s}}_{L^p_w(\R^n)}=\norm{T_s\circ \mathcal{R}_jv}_{L^p_w(\R^n)}\leq C\norm{v}_{L^p_w(\R^n)}<\infty,$$ so $\nabla^su\in L_w^p(\R^n;\R^n)$.
Now, we prove the other inclusion. Let $u\in \mathcal{S}(\R^n)\cap X_w^{s,p}(\R^n)$, which is enough since $C_c^\infty(\R^n)\subset \mathcal{S}(\R^n)\cap X_w^{s,p}(\R^n)$, and $C_c^\infty(\R^n)$ is dense in $X_w^{s,p}(\R^n)$
 by definition. We define $$v=G_s\circ \left(\operatorname{id}+\sum_{j=1}^n\mathcal{R}_j\frac{\partial^s}{\partial x_j^s}\right)u,$$ where $$\widehat{G_s\varphi}:=\frac{\left\langle\xi\right\rangle^s}{1+(2\pi|\xi|)^s}\widehat{\varphi},\,\varphi\in \mathcal{S}(\R^n).$$ Again, by Theorem \ref{Milhin} and the $L_w^p$-boundedness of the Riesz transform we have that $v\in L_w^p(\R^n)$, so it suffices to show that $u=\Lambda_s v$. Now, \begin{align*}
     \mathfrak{F}\{\Lambda_sv\}(\xi)&=(1+4\pi|\xi|^2)^{-s/2}\mathfrak{F}\{\lambda_s\}(\xi)\left(\mathfrak{F}\{u\}(\xi)+\sum_{j=1}^n\frac{-i\xi_j}{|\xi|}(2\pi)^si\xi_j|\xi|^{s-1}\mathfrak{F}\{u\}(\xi)\right)\\
     &=(1+4\pi|\xi|^2)^{-s/2}\frac{(1+4\pi|\xi|^2)^{s/2}}{1+(2\pi|\xi|)^s}\left(1+\frac{(2\pi)^s|\xi|^2}{|\xi|^{2-s}}\right)\mathfrak{F}\{u\}(\xi)\\
     &=\frac{1+(2\pi|\xi|)^2}{1+(2\pi|\xi|)^s}\mathfrak{F}\{u\}(\xi)=\mathfrak{F}\{u\}(\xi),
 \end{align*}
 which implies that $\Lambda_s v=u$, and thus the desired result follows.\qed\\
\section{Fractional weighted spaces on bounded domains}
\subsection{Interpolation on bounded domains and the fractional-weighted setting}
Since we want to establish the setting for the study of natural phenomena, we would like to adapt the definitions of the last section to suitable bounded domains. The development of weighted Bessel potential spaces on bounded domains was done by Frohlich in \cite{Frohlich2004} and Schumacher in \cite{Schumacher2009I, Schumacher2009II}. Before defining the Fractional-Weighted Sobolev spaces on bounded domains, we explore the interpolation nature of the weighted Bessel potential spaces in bounded domains.
\begin{Defi}\em
    A domain $\Omega\subset \R^n$ is said to be an \textit{extension domain} if for all $m\in \mathbb{N}$, there exists an extension operator $E_m$ such that for all $1<p<\infty$, $w\in A_p$ and $k=0,1,\ldots,m$, $E_m:W^{k,p}_w(\Omega)\to W^{k,p}_w(\R^n)$ is linear and bounded.
\end{Defi}
There are several classes of domains which are extension domains. We are mainly interested in Lipschitz domains, for whom existence of extension operators is guaranteed (see \cite[Theorem 2.3]{Frohlich2004} or \cite[Lemma 3.4]{Schumacher2009I}).\\

For an extension domain $\Omega$, we define the weighted Bessel potential space on $\Omega$ by $$H^{s,p}_w(\Omega):=\{u|_{\Omega}: u\in H^{s,p}_w(\R^n)\},\,s\in\R,\,1<p<\infty,\,w\in A_p,$$ with the norm $$\norm{u}_{H^{s,p}_w(\Omega)}:=\operatorname{inf}\{\norm{v}_{H^{s,p}_w(\R^n)}:v\in H^{s,p}_w(\R^n), v|_\Omega=u\}.$$
From now on, $\Omega$ will denote an open bounded subset of $\R^n$ with Lipschitz boundary. Since we have the existence of an extension operator for the weighted Sobolev spaces $W^{k,p}_w(\Omega)$, independent on the differentiability order $k$, we can easily derive the following equivalence interpolation equivalence. 
\begin{Teor}
    Let $s>0$, $1<p<\infty$, $w\in A_p$ and $k\in \Z$. Then, $$H^{s,p}_w(\Omega)=[L_w^p(\Omega),W^{k,p}_w(\Omega)]_{s/k},$$ with equivalence of the norms, where $k$ is the leas integer greater than $s$.
\end{Teor}
\noindent{}\textbf{Proof:} The proof is identical to the one of the unweighted case. See \cite[Theorem 2.22]{BellidoCuetoGarcia2025} for example. \qed

The definition of the Riesz fractional gradient requires the integral to be defined on the whole space, and there is no suitable version of it for bounded domains. Truncated versions of the fractional gradient have been introduced and developed in \cite{BellidoCuetoMoraCorral2023} by limiting the horizon of interaction between particles to a ball of radius $\delta>0$. However, this definition depends on the parameter $\delta$, so a canonical fractional gradient for $H^{s,p}(\Omega)$ is still missing, a fortiori to $H^{s,p}_w(\Omega)$. Hence, we will work on bounded domains with functions extended by zero to the whole space. We define those spaces as the completion of compactly supported smooth function on $\Omega$ under the $H^{s,p}_w(\R^n)$ norm. In particular, we introduce the spaces $$H^{s,p}_{0,w}(\Omega):=\overline{C_c^\infty(\Omega)}^{\norm{\cdot}_{H^{s,p}_w(\R^n)}},\,s\in (0,1),\,1<p<\infty,\,w\in A_p.$$ Those spaces are complete and reflexive closed subspaces of $H^{s,p}_w(\R^n)$, and clearly when $\Omega=\R^n$, $H^{s,p}_w(\R^n)=H^{s,p}_{0,w}(\R^n)$ by the density of $C_c^\infty(\R^n)$ on the first space. The relationship between $H^{s,p}_w(\Omega)$ and $H^{s,p}_{0,w}(\Omega)$ comes from their dual spaces. In particular we have the following result \cite[Lemma 3.3]{Schumacher2009I}:
\begin{Prop}\label{Duality}
    Let $s\in \R$, $1<p<\infty$ and $w\in A_p$. Then, we have $$\left(H^{s,p}_{0,w}(\Omega)\right)^*=H^{-s,q}_{w^*}(\Omega),\,1/p+1/q=1,\,w^*=w^{-1/(q-1)}.$$
\end{Prop}
\noindent{}\textbf{Proof:} Let $f\in H^{-s,q}_{w^*}(\Omega)$. Then, there exists $F\in H^{-s,q}_{w^*}(\R^n)$ such that the restriction in the sense of distribution of $F$ is $f$, i.e., $F|_{C^\infty_c(\Omega)}=f$. Now, \begin{align*}
    \norm{f}_{\left(H^{s,p}_{0,w}(\Omega)\right)^*}&=\operatorname{sup}_{\varphi\in C_c^\infty(\Omega): \norm{\varphi}_{H^{s,p}_{0,w}(\Omega)}=1}\left\langle f,\varphi\right\rangle\leq \operatorname{sup}_{\varphi\in \mathcal{S}(\R^n): \norm{\varphi}_{H^{s,p}_{w}(\R^n)}=1}\left\langle F,\varphi\right\rangle\\
    &=\norm{F}_{H^{-s,q}_{w^*}}\leq 2\norm{f}_{H^{-s,q}_{w^*}(\Omega)},
\end{align*} where we have used the fact that $\left(H^{s,p}_w(\R^n)\right)^*=H^{-s,q}_{w^*}(\R^n)$ isometrically. Hence, $f\in \left(H^{s,p}_{0,w}(\Omega)\right)^*$. For the other inclusion, given $g\in \left(H^{s,p}_{0,w}(\Omega)\right)^*$, we can extend it to an element $G\in \left(H^{s,p}_w(\R^n)\right)^*=H^{-s,q}_{w^*}(\R^n)$ by Hahn-Banach's theorem, with equality of the norms. Hence, an analogous computation yields that $$\norm{g}_{H^{-s,q}_{w^*}(\Omega)}\leq \norm{u}_{\left(H^{s,p}_{0,w}(\Omega)\right)^*},$$ and hence $g\in H^{-s,q}_{w^*}(\Omega)$.\qed\\

With this result, we are ready to prove that the spaces $H^{s,p}_{0,w}(\Omega)$ are interpolation spaces as well. 
\begin{Teor}\label{interpol}
    Let $s>0$, $1<p<\infty$ and $w\in A_p$. Then, $$H^{s,p}_{0,w}(\Omega)=\left[L^p_w(\Omega),W^{k,p}_{0,w}(\Omega)\right]_{s/k},$$ where $k$ is the least integer greater than $s$.
\end{Teor}
\noindent{}\textbf{Proof:}
By the reiteration theorem for the complex method we have that \begin{align*}
    H^{-s,q}_{w^*}(\Omega)=\left[L^{q}_{w^*}(\Omega),W^{-k,q}_{w^*}(\Omega)\right]_{s/k},
\end{align*}
and by the duality theorem $$\left[L^{q}_{w^*}(\Omega),W^{-k,q}_{w^*}(\Omega)\right]_{s/k}=\left[\left(L^{p}_{w}(\Omega)\right)^*,\left(W^{k,p}_{w}(\Omega)\right)^*\right]_{s/k}=\left(\left[L^{p}_{w}(\Omega),W^{k,p}_{0,w}(\Omega)\right]_{s/k}\right)^*,$$ hence $$\left(\left[L^{p}_{w}(\Omega),W^{k,p}_{0,w}(\Omega)\right]_{s/k}\right)^*=H^{-s,q}_{w^*}(\Omega),$$ thus $$H^{s,p}_{0,w}(\Omega)=\left[L^{p}_{w}(\Omega),W^{k,p}_{0,w}(\Omega)\right]_{s/k},$$ as we wanted to prove. \qed

With $H^{s,p}_{w,0}(\Omega)$ on mind we introduce the following spaces by means of the Riesz fractional gradient.
\begin{Defi}\em[Fractional-Weighted Sobolev spaces]
Let $s\in (0,1)$, $1<p<\infty$ and $w\in A_w$. We define the space $X^{s,p}_{0,w}(\Omega)$ as the completion of $C_c^\infty(\Omega)$ under the norm $\norm{\cdot}_{\nabla^s,p,w}$, i.e., $$X^{s,p}_{0,w}(\Omega):=\overline{C_c^\infty(\Omega)}^{\norm{\cdot}_{\nabla^s,p,w}}.$$
\end{Defi}
For an element $u\in C_c^\infty(\Omega)$, by $\norm{u}_{\nabla^s,p,w}=\norm{u}_{L^p_w(\R^n)}+\norm{\nabla^su}_{L^p_w(\R^n)}$, we mean the norm of the extension by zero of $u$ to the whole $\R^n$.\\
As we stated for the case $\Omega=\R^n$, to be consequent with our definition, we have to extend the definition of $\nabla^su$ for functions $u\in X_{0,w}^{s,p}(\Omega)$ by continuity. Since by definition, for every $u\in X_{0,w}^{s,p}(\Omega)$ there exists a sequence $\{u_m\}\subset C_c^\infty(\Omega)$ such that $u_m\to u$ in $L^p_w(\R^n)$ and $\{\nabla^su_m\}$ is Cauchy in $L^p_w(\R^n),\R^n)$, we define $$\nabla^su:=L^p_w(\R^n;\R^n)-\lim_{m\to \infty} \nabla^su_m.$$ This characterization also allows us to extend the integration by parts formula for functions in $X_{0,w}^{s,p}(\Omega)$.\\

Note that since we have the followings equivalences of the norms $$\norm{\cdot}_{\nabla^s,p,w}\sim \norm{\cdot}_{H^{s,p}_w(\R^n)}\sim \norm{\cdot}_{[L^p_w(\R^n),W^{1,p}_w(\R^n)]_s},$$ we have that $X_{0,w}^{s,p}(\Omega)$ is equivalent to $H^{s,p}_{0,w}(\R^n)$. Moreover, $X_{0,w}^{s,p}(\Omega)$ is equivalent to the complex interpolation of the couple $\left(L^p_w(\Omega),W^{1,p}_{0,w}(\Omega)\right)$. However, we will maintain the notation $X_{0,w}^{s,p}(\Omega)$ instead of $H^{s,p}_{0,w}(\Omega)$ to remark that we are working with the fractional gradient. 
\subsection{Continuous embeddings and interpolation inequalities}
Once we have established the function space structure for our spaces, we study some important embeddings and inequalities. We first prove the analogous to the classical Sobolev inequality for the fractional-weighted setting. It follows from the following weighted estimate for the Riesz potential \cite[Theorem 2.2.1.]{turesson2000}. We first introduce the following notations: Let $0<s<n$ and $1<p<\infty$ such that $sp<n$. We define the \textit{fractional Sobolev conjugate} $p_s^*$ as $$p_s^*=\frac{np}{n-sp}.$$
In the unweighted case, the Sobolev embedding on the subcritical case says that \cite[Theorem 3.21]{BellidoGarcia2025} $$H^{s,p}(\R^n)\xhookrightarrow{}L^{p_s^*}(\R^n),$$ and hence $$H^{s,p}_0(\Omega)\xhookrightarrow{}L^{p_s^*}(\Omega).$$ Moreover, by the Hardy-Littlewood-Sobolev inequality, that says that $I_s:L^p(\R^n)\to L^{p_s^*}(\R^n)$ continuously, for every $u\in H^{s,p}_0(\Omega)$, $$\norm{u}_{L^{p_s^*}(\Omega)}=\norm{I_s\mathcal{R}\cdot D^su}_{L^{p_s^*}(\R^n;\R^n)}\leq C\norm{D^su}_{L^p(\R^n;\R^n)}.$$ The weighted case is quite different since the Riesz potential also affects to the weight $w$. For $w\in A_{p_s^*/p'+1}$, where $1/p+
1/p'=1$. We denote $$w_{s,p}:=w^{(n-sp)/n}=w^{p/p_s^*}.$$ The weighted estimate for the Riesz potential now reads as:
\begin{Teor}\label{WeigthedRiesz}
    Let $0<s<n$ and $1<p<\infty$ such that $sp<n$. Let $w\in A_{p_s^*/p'+1}$. Then, $$\left(\int_{\R^n}|I_su|^{p_s^*}w\,dx\right)^{1/p_s^*}\leq C\left(\int_{\R^n}|f|^pw_{s,p}\,dx\right)^{1/p},\,u\in L^p_{w_{s,p}}(\R^n),$$ where $C$ only depends on $s,n,p$ and the $A_p$-constant of $w$.
\end{Teor}
From here we can establish the following Sobolev embedding: 
\begin{Teor}[Sobolev inequality]\label{EmbeddingI}
    Let $s\in (0,1)$ and $1<p<\infty$ such that $sp<n$. Let $w\in A_{p_s^*/p'+1}$. Then, for every $u\in X_{0,{w_{s,p}}}^{s,p}(\Omega)$,$$\norm{u}_{L_w^{p_s^*}(\Omega)}\leq C\norm{\nabla^su}_{L_{w_{s,p}}^p(\R^n)}.$$
\end{Teor}
\noindent{}\textbf{Proof:} Let $u\in C_c^\infty(\Omega)$. From Theorem \ref{WeigthedRiesz}, $$\norm{I_sv}_{L_w^{p_s^*}(\R^n)}\leq C\norm{v}_{L_{w_{s,p}}^p(\R^n)},$$ for all $v\in L_{w_{s,p}}^p(\R^n)$.
By the fractional fundamental theorem of calculus stated in Proposition \ref{DsProps}, $u=I_s\left(\mathcal{R}\cdot \nabla^su\right)$. By the weighted $L^p$-boundedness of the Riesz transform, choosing $v=\mathcal{R}\cdot \nabla^su$, we get that $$\norm{u}_{L^{p_s^*}_w(\Omega)}=\norm{u}_{L^{p_s^*}_w(\R^n)}=\norm{I_sv}_{L^{p_s^*}_w(\R^n)}\leq C\norm{\mathcal{R}\cdot \nabla^su}_{L^p_{w_{s,p}}(\R^n)}\leq C\norm{\nabla^su}_{L^p_{w_{s,p}}(\R^n)}.$$ Now, the result follows from density. \qed

As a corollary, we obtain the following continuous embedding which generalizes \cite[Theorem 3.21 (1)]{BellidoGarcia2025}:\begin{Coro}
    Let $s\in (0,1)$, $1<p<\infty$ such that $sp<n$ and $w\in A_{p_s^*/p+1}$. Then, we have that $$X_{0,w_{s,p}}^{s,p}(\Omega)\xhookrightarrow{}L^{p_s^*}_w(\Omega).$$
\end{Coro}
If we allow the weight $w$ to change in the target space of $I_s$ acting on $L^p_w$, a more general inequality for the Riesz potential can be obtained. In particular, we have the following result \cite[Theorem 1(B)]{SawyerWheeden1992}:
\begin{Teor}\label{RieszII}
    Let $0<s<n$ and $1<p<q<\infty$. Let $w\in A_p$ and $v\in A_q$ such that $$|Q|^{s/n-1}\left(\int_Qv(x)\,dx\right)^{1/q}\left(\int_Qw^*(x)\,dx\right)^{1/p'}<C,$$ for every cube $Q\subset \R^n$, with $C>0$ independent of the cube. Then, for every $u\in L_w^p(\R^n)$, $$\norm{I_su}_{L^q_v(\R^n)}\leq C\norm{f}_{L^p_w(\R^n)}.$$
\end{Teor}
Again, by means of the fractional theorem of calculus we get the following:
\begin{Teor}\label{EmbeddingII}
    Let $0<s<n$ and $1<p<q<\infty$. Let $w\in A_p$ and $v\in A_q$ such that $$|Q|^{s/n-1}\left(\int_Qv(x)\,dx\right)^{1/q}\left(\int_Qw^*(x)\,dx\right)^{1/p'}<C,$$ for every cube $Q\subset \R^n$, with $C>0$ independent of the cube. Then, we have the embedding $$X_{0,w}^{s,p}(\Omega)\xhookrightarrow{}L^q_v(\Omega),$$ with the inequality $$\norm{u}_{L^p_w(\Omega)}\leq C\norm{\nabla^su}_{L^q_v(\Omega)}.$$
\end{Teor}

In the case that we impose $w=v\in A_p\subset A_q$, the condition 
$$|Q|^{s/n-1}\left(\int_Qw(x)\,dx\right)^{1/q}\left(\int_Qw^*(x)\,dx\right)^{1/p'}<C,$$ is equivalent to having $$|Q|^{s/n}w(Q)^{1/q-1/p}\leq C,$$ as it was proved in \cite{Frohlich2004}  Also those conditions do not have to verify for every cube $Q\subset \R^n$. In fact, it is enough to have it only for cubes $Q$ contained in a neighborhood of $\overline{\Omega}$ whenever it is an extension domain, as Lipschtiz ones (see \cite[pp. 156]{Frohlich2004}).
Hence, as a corollary of Theorem \ref{EmbeddingII} we obtain the following:
\begin{Coro}
    Let $s\in (0,1)$, $1<p<q<\infty$ and $w\in A_p$ such that $$|Q|^{s/n}w(Q)^{1/q-1/p}\leq C,$$ for every cube $Q\subset U\subset \overline{\Omega}$. Then, $$X_{0,w}^{s,p}(\Omega)\xhookrightarrow{}L^q_w(\Omega).$$
\end{Coro}

Now, ff we choose $1<p\leq q\leq r<\infty$, $w\in A_p$ and $s\in (0,1)$, with $$\frac{1}{r}\geq \frac{1}{q}-\frac{s}{np},$$ it is proven in \cite[Lemma 5.3]{Schumacher2009I} that this implies $$|Q|^{s/n}w(Q)^{1/r-1/q}\leq C,$$ for all cubes $Q\subset U$ with $\overline{\Omega}\subset U$. Hence, we obtain another direct corollary for Theorem \ref{RieszII}:
\begin{Coro}
    Let $1<p\leq q\leq r<\infty$ and $s\in (0,1)$ such that $$\frac{1}{r}\geq \frac{1}{q}-\frac{s}{np}.$$ Then, for every $w\in A_p$, $$X_{0,w}^{s,q}(\Omega)\xhookrightarrow{}L_w^r(\Omega).$$ 
\end{Coro}

For $1<p<q<\infty$, we say that a weight $w$ satisfies the $A_{p,q}$-condition, denoted as $w\in A_{p,q}$, if $$[w]_{p,q}:=\operatorname{sup}_Q\left(\frac{1}{|Q|}\int_Q w(x)\,dx\right)\left(\frac{1}{|Q|}\int_Q w(x)^{-p'/q}\,dx\right)^{q/p'}<\infty,$$ which is related with the restrictions over the weights for the latest results. We want to note that we imposed this restriction over the weight in Theorem \ref{EmbeddingI}, since $w\in A_{p,q}$ if and only if $w\in A_r$, $r=1+q/p'$ . Hence, Theorem \ref{EmbeddingI} could be re stated as follows:
\begin{Teor}
    Let $s\in (0,1)$ and $1<p<\infty$ such that $sp<n$. Then, for every $q\in [1,p_s^*]$, $$X^{s,p}_{0,w_{s,p}}(\Omega)\xhookrightarrow{}L^{q}_w(\Omega),\,w\in A_{p,p_s^*}.$$
\end{Teor}

As in \cite[Theorem 2.9]{BellidoCuetoMoraCorral2021} for $H_0^{s,p}(\Omega)=X_{0,1}^{s,p}(\Omega)$ and \cite[Theorem 2]{campos2024} for fractional Orlicz-Sobolev spaces, we have the following Poincar\'e inequality:
\begin{Teor}[Poincar\'e inequality]\label{Poincare} Let $\Omega\subset \R^n$ a bounded open set, $s\in (0,1)$, $1<p<\infty$ and $w\in A_p$. Then, for every $u\in X_{0,w}^{s,p}(\Omega)$, there exists a constant depending on $\Omega$, $n$, and $[w]_p$ such that  $$\norm{u}_{L^p_w(\Omega)}\leq \frac{C}{1-2^{-s}}\norm{\nabla^su}_{L^p_w(\R^n;\R^n)}.$$
\end{Teor}
\noindent{}\textbf{Proof:}
We follow the strategy of \cite[Theorem 2]{campos2024}. By density, it is enough to prove the result for $u\in C_c^\infty(\Omega)$. Consider $r>0$ large enough to have $\Omega\subset B_r(0)$. By the fractional fundamental theorem of calculus, \begin{align*}u(x)&=c_{n,-s}\int_{\R^n}\nabla^su(y)\cdot \frac{x-y}{|x-y|^{n-s+1}}\,dy\\
&=c_{n,-s}\left(\int_{B_{2r}}\nabla^su(y)\cdot \frac{x-y}{|x-y|^{n-s+1}}\,dy+\int_{B_{2r}^c}\nabla^su(y)\cdot \frac{x-y}{|x-y|^{n-s+1}}\,dy\right),\end{align*} and hence $$\norm{u}_{L^p_w(\Omega)}\leq c_{n,-s}\left(\norm{\int_{B_{2r}}\frac{|\nabla^su(y)|}{|x-y|^{n-s}}\,dy}_{L^p_w(\Omega)}+\norm{\int_{B_{2r}^c}\frac{|\nabla^su(y)|}{|x-y|^{n-s}}\,dy}_{L^p_w(\Omega)}\right).$$
We now proceed to estimate each integral separately. First, let $y\in B_{2r}^c$. By definition, $$|\nabla^su(y)|\leq c_{n,s}\int_{\Omega}\frac{|u(z)|}{|y-z|^{n+s}}\,dz.$$ Since $z\in \Omega\subset B_r$, and $y\in B_{2r}^c$, $|y-z|\geq \frac{|y|}{2}$, hence. On the other hand, by H\"older's inequality, 
\begin{align*}
    \int_{\Omega}|u(z)|\,dz=\int_\Omega |u(z)|w^{1/p}(z)w^{-1/p}(z)\,dz\leq \left(\int_\Omega |u(z)|^pw(z)\,dz\right)^{1/p}\left(\int_\Omega w(z)^{-q/p}\,dz\right)^{1/q},
\end{align*}
with $1/q+1/p=1$. Observe that $-q/p=-1/(p-1)$, and hence 
\begin{align*}
    |\nabla^su(y)|\leq \frac{2^{n+s}c_{n,s}}{|y|^{n+s}}\int_\Omega|u(z)|\,dz\leq \frac{2^{n+s}c_{n,s}}{|y|^{n+s}}\left(\int_{\Omega}|u(z)|^pw(z)\,dz\right)^{1/p}\left(\int_\Omega w^*(z)\,dz\right)^{1/q},
\end{align*}
so $$|\nabla^su(y)|\leq \frac{2^{n+s}c_{n,s}}{|y|^{n+s}}\norm{u}_{L^p_w(\Omega)}w^*(\Omega)^{1/q}.$$ Hence, \begin{align*}\int_{B_{2r}^c}\frac{|\nabla^su(y)|}{|x-y|^{n-s}}\,dy&\leq \frac{2^{n+s}c_{n,s}}{|y|^{n+s}}\norm{u}_{L^p_w(\Omega)}w^*(\Omega)^{1/q}\int_{B_{2r}^c}\frac{1}{|y|^{n+s}|x-y|^{n-s}}\,dy\\
&\leq \frac{2^{3-n}c_{n,s}\omega_nr^{-n}}{n}w^*(\Omega)^{1/q}\norm{u}_{L^p_w(\Omega)},\end{align*} thus \begin{align*}\norm{\int_{B_{2r}^c}\frac{|\nabla^su(y)|}{|x-y|^{n-s}}\,dy}_{L^p_w(\Omega)}&\leq \frac{2^{3-n}c_{n,s}\omega_nr^{-n}}{n}w^*(\Omega)^{1/q}\norm{u}_{L^p_w(\Omega)}\left(\int_\Omega w(x)\,dx\right)^{1/p}\\
&=\frac{2^{3-n}c_{n,s}\omega_nr^{-n}}{n}w^*(\Omega)^{1/q}\norm{u}_{L^p_w(\Omega)}w(\Omega)^{1/p}.\end{align*}
Observe that $$w(\Omega)^{1/p}w^*(\Omega)^{1/q}\leq |\Omega|[w]_p^{1/p},$$ and hence $$\norm{\int_{B_{2r}^c}\frac{|\nabla^su(y)|}{|x-y|^{n-s}}\,dy}_{L^p_w(\Omega)}\leq \frac{2^{3-n}c_{n,s}\omega_nr^{-n}}{n}|\Omega|[w]_p^{1/p}\norm{u}_{L^p_w(\Omega)}.$$
For the other integral we use \cite[Lemma 6.1.4]{DieningEtAl2011}, which says that for every $x\in \R^n$, $\delta>0$, $0<s<n$, and $f\in L_{loc}^1(\R^n)$, we have that $$\int_{B(x,\delta)}\frac{|f(y)|}{|x-y|^{n-s}}\,dy\leq C\frac{\delta^s}{1-2^{-s}}Mf(x),$$where $C$ does not depend on $s$,
 and hence \begin{align*}
     \int_{B_{2r}}\frac{|\nabla^su(y)|}{|x-y|^{n-s}}\,dy\leq \int_{|x-y|\leq 3r}\frac{|\nabla^su(y)|}{|x-y|^{n-s}}\,dy\leq C\frac{(3r)^s}{1-2^{-s}}M\nabla^su(x),
 \end{align*} and since $M:L_w^p(\Omega)\to L_w^p(\Omega)$ continuously, \begin{align*}
     \norm{ \int_{B_{2r}}\frac{|\nabla^su(y)|}{|x-y|^{n-s}}\,dy}_{L^p_w(\Omega)}\leq C'3r\frac{1}{1-2^{-s}}\norm{\nabla^su}_{L^p_w(\Omega)},
 \end{align*}
where $C'$ does not depend on $s$. Combining both estimates and using the fact that $c_{n,s}$ is uniformly bounded on the parameter $s\in [-1,1]$ (see \cite[Lemma 2.4]{BellidoCuetoMoraCorral2021}), we obtain that there exists a positive constant $C''$, not depending on $s$, such that $$\norm{u}_{L^p_w(\Omega)}\leq C''\left(r^{-n}\norm{u}_{L^p_w(\Omega)}+\frac{r}{1-2^{-s}}\norm{\nabla^su}_{L^p_w(\Omega)}\right).$$ The result now follows choosing $e$ large enough to have $C''r^{-n}\leq 1/2$. \qed\\

Note that the Poincar\'e inequality allows us to characterize the elements in $\left(X^{s,p}_{0,w}(\Omega)\right)^*=H^{-s,q}_{w^*}(\Omega)$. In fact, for every $F\in H^{-s,q}_{w^*}(\Omega)$, there exists $f\in L^{q}_{w^*}(\Omega)$ such that $$Fu=\int_{\R^n}f\cdot \nabla^su,\,u\in H^{s,p}_{0,w}(\Omega).$$

Let $0\leq r\leq s\leq t\leq 1$, with $r\not =t$.
Observe that the interpolation identity $$\left[H^{r,p}_{0,w}(\Omega),H^{t,p}_{0,w}(\Omega)\right]_\theta=H^{s,p}_{0,w}(\Omega),\,\theta=\frac{s-r}{t-r},$$ which follows from Theorem \ref{ComplexProps} (8), implies the following  Gagliardo-Nirenberg inequality for functions $u\in H_{0,w}^{t,p}(\Omega)$,  $$\norm{u}_{H_{0,w}^{s,p}(\Omega)}\leq C \norm{u}_{H_{0,w}^{r,p}(\Omega)}^{1-\theta}\norm{u}_{H_{0,w}^{t,p}(\Omega)}^\theta,\,\theta=\frac{s-r}{t-r},$$ which follows as well from Theorem \ref{ComplexProps} (5). In view of this we would expect some interpolation inequality regarding the $L^p_w-$norm of the Riesz fractional gradient. We present the following result generalizing 
\cite[Theorem 13]{Brue2022}.
\begin{Teor}
    Let $0\leq r\leq s\leq t<1$, $1<p<\infty$ and $w\in A_p$. Then, for every $u\in H_{0,w}^{t,p}(\Omega)$, $$\norm{\nabla^su}_{H_{0,w}^{s,p}(\Omega)}\leq C\norm{\nabla^ru}_{H_{0,w}^{r,p}(\Omega)}^{1-\theta}\norm{\nabla^tu}_{H_{0,w}^{t,p}(\Omega)}^\theta,\,\theta=\frac{s-r}{t-r},$$ where $C$ does not depend on $r,s,t$.
\end{Teor}
\noindent{}\textbf{Proof:}
Let us define the function $$m_{t,s}(\xi):=\frac{(2\pi|\xi|)^s}{1+(2\pi|\xi|)^t},\xi\in \R^n.$$ By the computations done in \cite[Lemma 3]{campos2024}, is easy to see that $$\sup_{\alpha\in \N_0:|\alpha|\leq \left[\frac{n}{2}\right]+1}\sup_{\xi\in \R^n\setminus\{0\}}\left|\xi^{|\alpha|}\partial^\alpha_\xi m_{t,s}(\xi)\right|<\infty.$$ Hence, Theorem \ref{Milhin} implies that the operator $$\widehat{T\varphi}:=m_{s,t}\widehat{\varphi},\,\varphi\in \mathcal{S}(\R^n),$$ extends to a bounded operator from $L_w^p(\R^n)\to L_w^p(\R^n)$. Moreover, note that for every $u\in C_c^\infty(\Omega)$, $$\widehat{(-\Delta)^{s/2}u}=(2\pi |\xi|)^s\widehat{u}=\frac{(2\pi |\xi|)^s}{1+(2\pi|\xi|)^t}\left(1+(2\pi|\xi|)^t\right)\widehat{u}=\mathfrak{F}\Bigg\{T_{t,s}\circ \left(id+(-\Delta)^{t/2}\right)u\Bigg\},$$ and hence $(-\Delta)^{s/2}u=T_{t,s}\circ \left(id+(-\Delta)^{t/2}\right)u$. Now, observe that from the fractional fundamental theorem of calculus, $$(-\Delta)^{\alpha/2}\nabla^su=(-\Delta)^{(\alpha+s)/2}\mathcal{R}u,\,\alpha\in (0,1).$$ Hence, we can estimate \begin{align*}\norm{(-\Delta)^{s/2}u}_{L^p_w(\R^n)}&=\norm{T_{t,s}\circ \left(id+(-\Delta)^{t/2}\right)u}_{L^p_w(\R^n)}\leq C\norm{\left(id+(-\Delta)^{t/2}\right)u}_{L^p_w(\R^n)}\\
&\leq C \norm{u}_{L^p_w(\R^n)}+\norm{(-\Delta)^{t/2}u}_{L^p_w(\R^n)},\end{align*} where the constant $C$ does not depend on $t$ and $s$. Hence, $$\norm{\nabla^su}_{L^p_w(\R^n)}\leq C \norm{\mathcal{R}u}_{L^p_w(\R^n)}+\norm{\nabla^tu}_{L^p_w(\R^n)}\leq C \norm{u}_{L^p_w(\R^n)}+\norm{\nabla^tu}_{L^p_w(\R^n)},$$and performing a dilation and optimizing the last expression we find that 
$$\norm{\nabla^su}_{L_w^p(\R^n;\R^n)}\leq C\norm{u}_{L_w^p(\R^n;\R^n)}^{\frac{t-s}{t}}\norm{\nabla^tu}_{L_w^p(\R^n;\R^n)}^{\frac{s}{t}},$$ which proves the case $r=0$. Now, since $$(-\Delta)^{(s-r)/2}\nabla^ru=\mathcal{R}(-\Delta)^{s/2}u=\nabla^su,$$ in an analogous way we can prove the case $r>0$, \begin{align*}
    \norm{\nabla^su}_{L^p_w(\R^n)}&=\norm{(-\Delta)^{(s-r)/2}\nabla^ru}_{L^p_w(\R^n)}\leq C \norm{\nabla^ru}_{L^p_w(\R^n)}+\norm{(-\Delta)^{(t-r)/2}\nabla^ru}_{L^p_w(\R^n)}\\
&\leq C\norm{\nabla^ru}_{L^p_w(\R^n)}+C\norm{\nabla^tu}_{L^p_w(\R^n)},\end{align*} and hence 
$$\norm{\nabla^su}_{L^p_w(\R^n;R^n)}\leq C \norm{\nabla^ru}_{L^p_w(\R^n)}^{1-\theta}\norm{\nabla^tu}_{L^p_w(\R^n)}^\theta,$$ for $\theta=\frac{s-r}{t-r}$. Finally, we can extend the result by density for functions $u\in H^{s,p}_{0,w}(\Omega)$. \qed

In \cite{Duarte2023}, similar Gagliardo-Nirenberg inequalities are obtained for different fractional operators, in particular for the Bessel potential and the fractional Laplacian. 
\subsection{Compactness results}
To conclude the study of basic results for our spaces $X_{0,w}^{s,p}(\Omega)$, we prove the compactness of the embedding $$X_{0,w}^{s,p}(\Omega\xhookrightarrow{}L^p_w(\Omega).$$ We provide a straightforward proof based on the compact embedding for the local case $s=1$ and the complex interpolation method.
\begin{Teor}
Let $\Omega\subset \R^n$ a Lipschitz domain, $s\in (0,1)$, $1<p<\infty$ and $w\in A_p$. Then, the embedding $$X_{0,w}^{s,p}(\Omega)\xhookrightarrow{}L^p_w(\Omega),$$
    is compact.
\end{Teor}
\noindent{}\textbf{Proof:}
By \cite[Theorem 2.3]{Frohlich2007}, the embedding $$W_{0,w}^{1,p}(\Omega)\xhookrightarrow{}L^p_w(\Omega),$$ is compact. On the other hand, by Theorem \ref{interpol}, $$X_{0,w}^{s,p}(\Omega)=\left[L^p_w(\Omega),W_{0,w}^{1,p}(\Omega)\right]_s,$$ and hence interpolating the embeddings $L^p_w(\Omega)\xhookrightarrow{}L^p_w(\Omega)$ and $W^{1,p}_{0,w}(\Omega)\xhookrightarrow{}L^p_w(\Omega)$, we get that $$X_{0,w}^{s,p}(\Omega)\xhookrightarrow{}L^p_w(\Omega).$$ The key is that given two compatible couples of Banach spaces $(E_0,E_1)$ and $(F_0,F_1)$ such that $E_i\xhookrightarrow{}F_i$, $i=0,1$, with $E_1\xhookrightarrow{}F_1$ being compact, the complex interpolation method preserve the compactness if $E_1$ has the UMD condition (see \cite{CwikelKalton1995}), i.e., in that case the embedding $$[E_0,E_1]_\theta\xhookrightarrow{}[F_0,F_1]_\theta,\,\theta\in(0,1),$$ is compact. Since $W_{0,w}^{1,p}(\Omega)$ has the UMD condition, we conclude that $$X_{0,w}^{s,p}(\Omega)\xhookrightarrow{}L^p_w(\Omega),$$ compactly.\qed

If less regularity for the domain $\Omega$ is assumed, and hence we do not have the equivalence of our space with the complex interpolation one, the usual path to prove the compact embedding would be adapting the ideas of \cite[Theorem 2.3]{Frohlich2007} with the integral estimates from \cite[Theorem 2.1]{ShiehSpector2018}. Another possibility would be by means of the weighted Fr\'echet-Kolmogorov theorem from \cite[Theorem 5]{ClopCruz2013} and an estimate on the translations for the $L^p_w$-norm analogous to the one obtained in \cite[Proposition 3.3]{BellidoCuetoGarcia2025} for the unweighted case. 
\section{A Fractional nonlinear Degenerate Elliptic problem}
In this section we will address the problem of existence and uniqueness of weak solutions of the following class of degenerate elliptic PDEs on the fractional setting. Let $\Omega\subset \R^n$ an open bounded domain, $s\in (0,1)$ and $1<p<\infty$. We consider the following nonlinear problem: 
\begin{equation}\label{eqn:P}\begin{cases}
-\operatorname{div}^s\left(w(x)|\nabla^su|^{p-2}\nabla^s u\right)&=f,\,x\in \Omega,\\
    \hfill u&=0,\,x\in \Omega^c.
\end{cases}\end{equation}
Here, $f$ is some functional and  
$w\in A_p$. \\
The unweighted case, i.e., $w(x)=1$, the operator $$H^{s,p}_0(\Omega)\to H^{-s,p'}(\Omega); u\mapsto-\operatorname{div}^s(|\nabla^su|^{p-2}\nabla^su),$$ is the so-called $H^{s,p}$-Laplacian introduced in \cite{Schikorra2018}. This operator arises naturally as the first variation of energies of the norm $\norm{\nabla^su}_p$, and seems as the suitable fractional generalization of the classical $p$-Laplacian $$-\operatorname{div}(|\nabla u|^{p-2}\nabla u),$$ instead of the so called \textit{fractional $p$-Laplacian}, which arises form the first variation of energies of the semi-norm of the \textit{Sobolev-Slobodeckij} spaces $W^{s,p}$. It is worth mentioning that both objects coincide only in the Hilbertian case $p=2$ (see \cite[Section 3.3]{BellidoGarcia2025} and \cite[Section 2.3.3]{Campos2021}). 

Suppose that $u$ solves \ref{eqn:P} and let $v\in C_c^\infty(\Omega)$. Multiplying both sides in \ref{eqn:P}, integrating over the whole space and formally using the integration by parts formula we get that $u$ is a weak solution of the problem if it satisfies $$\int_{\R^n}w(x)|\nabla^su|^{p-2}\nabla^su\cdot \nabla^sv\,dx=\int_\Omega fv\,dx,\,\forall v\in C_c^\infty (\Omega).$$ 
Since $w$ is a weight in the class $A_p$, it is natural to search for weak solutions in the weighted space $X^{s,p}_{0,w}(\Omega)$ and hence, in order to the expression $\int_\Omega fv$ to make sense we will assume that $f\in \left(X^{s,p}_{0,w}(\Omega)\right)^*=H^{-s,p'}_{w^*}(\Omega)$, with $w^*=w^{-1/(p-1)}\in A_{p'}$, $1/p+1/p'=1$. Hence, for $v\in X^{s,p}_{0,w}(\Omega)$, we have that \begin{align*}
    \int_\Omega fv\,dx=\int_\Omega fw(x)^{-1/p}w(x)^{1/p}v\,dx\leq \left(\int_\Omega |f|^{p'}w(x)^{-p'/p}\,dx\right)^{1/p'}\left(\int_\Omega |v|^pw(x)\,dx\right)^{1/p},
\end{align*} and since $$w^{-p'/p}=w^{-1/(p'-1)}=w^*,$$ we have that $$\int_\Omega fv\,dx\leq \norm{f}_{L^{p'}_{w^*}(\Omega)}\norm{v}_{L^p_w(\Omega)}<\infty.$$
\begin{Teor}\label{P}
    Let $\Omega\subset \R^n$ an open bounded domain $s\in (0,1)$ and $1<p<\infty.$ Let $w\in A_p$ and $f\in H^{-s,p'}_{w^*}(\Omega)$. Then, there exists an unique weak solution $u\in X_{0,w}^{s,p}(\Omega)$ for the problem \ref{eqn:P}, i.e., there exists an unique $u\in X_{0,w}^{s,p}(\Omega)$ such that $$\int_{\R^n}w(x)|\nabla^su|^{p-2}\nabla^su\cdot \nabla^sv\,dx=\int_\Omega fv\,dx,\,\forall v\in X^{s,p}_{0,w}(\Omega).$$
\end{Teor}
In order to prove the theorem we have to recall some aspects of the theory of monotone operators and its link with existence of minimizers for variational problems.

Let $X$ some Banach space and $X^*$ its dual. An operator $T:X\to X^*$ is called \textit{monotone} if $$\langle Tx-Ty,x-y\rangle\geq 0,\,\forall x,y\in X,$$ where $\langle,\rangle$ denotes the natural pairing in $X^*$. If for every $x\not =y$, 
$$\langle Tx-Ty,x-y\rangle> 0,$$ the operator is called \textit{strictly monotone}. The operator will be called \textit{coercive} if $$\lim_{\norm{x}_X}\frac{
\langle Tx,x\rangle}{\norm{x}_X}=\infty.$$  
Note that this definition of coercivity does not conflict with the one typically used for problems in Hilbert spaces. In fact, if we want to apply Lax-Milgram lemma to a bilinear form $a:H\times H\to \R$, whee $H$ is a Hilbert space, we require the form $a$ to be coercive in the sense that $$a(u,u)\geq c\norm{u}_H^2,$$ which implies, $$\lim_{\norm{x}_X\to \infty}\frac{a(u,u)}{\norm{u}_H}\geq c\lim_{\norm{u}_H\to \infty}\norm{u}_H=\infty.$$ Since the bilinear forms $a$ usually take the form of $\langle Tu,v\rangle$ for some $T:H\to H^*$, both definitions are the same.

We say that $T$ is \textit{demicontinuous} if for every $(x_n)$ converging to $x\in X$, we have that $\langle Tx_n,y\rangle \to \langle Tx,y\rangle$, for every $y\in X$, as $n\to \infty$. For those class of operators acting on separable and reflexive spaces, we have the following theorem , which states an existence result for problems on non Hilbert spaces:
\begin{Teor}[Browder-Minty]\label{BW}
    Let $X$ a separable, reflexive Banach space, and let $T:X\to X^*$ be a bounded, monotone, coercive and demicontinuous operator. Then, for all $f\in X^*$, there exists a solution $u\in X$ to the problem $$Tx=f.$$ If $T$ is strictly monotone, the solution is unique.
\end{Teor}
\noindent{}\textbf{Proof:} See \cite[Theorem 9.14]{Ciarlet2013} for a proof.

\noindent{}\textbf{Proof of Theorem \ref{P}:} Let the operator $$(-\Delta)_{p,w}^s:X^{s,p}_{0,w}(\Omega)\to H^{-s,p'}_{w^*}(\Omega),$$ defined as $$(-\Delta)_{p,w}^su:=-\operatorname{div}^s\left(w(x)|\nabla^su|^{p-2}\nabla^su\right).$$ The pairing is defined in a natural way by the integration by parts formula as $$\langle(-\Delta)_{p,w}^su,v\rangle=\int_{\R^n}w(x)|\nabla^su|^{p-2}\nabla^su\cdot \nabla^sv\,dx,$$ for all $v\in X^{s,p}_{0,w}(\Omega)$. Note that the pairing is well defined since by H\"older's inequality \begin{align*}
    |\langle(-\Delta)_{p,w}^su,v\rangle|&\leq \int_{\R^n}w(x)|\nabla^su|^{p-1}|\nabla^sv|\,dx=\int_{\R^n}w(x)^{\frac{p-1}{p}}|\nabla^su|^{p-1}w(x)^{1/p}|\nabla^sv|\,dx\\
    &\leq \left(\int_{\R^n}|\nabla^su|^pw(x)\,dx\right)^{\frac{p-1}{p}}\left(\int_{\R^n}|\nabla^s v|^pw(x)\,dx\right)^{1/p}=\norm{\nabla^s u}_{L^p_w}^{p-1}\norm{\nabla^sv}_{L^p_w}.
\end{align*}
The coercivity follows from the weighted Poincar\'e inequality, since $$\langle(-\Delta)_{p,w}^su,u\rangle=\int_{\R^n}w(x)|\nabla^su|^p\,dx=\norm{\nabla^su}_{L^p_w(\R^n;\R^n)}^p\geq c\norm{u}_{X_{0,w}^{s,p}(\Omega)}^p,$$ and hence $$\lim_{\norm{u}_{X_{0,w}^{s,p}(\Omega)}\to \infty}\frac{\langle(-\Delta)_{p,w}^su,u\rangle}{\norm{u}_{X_{0,w}^{s,p}(\Omega)}}=\infty,$$ since $p>1$.
For the monotonicity, let $u,v\in X_{0,w}^{s,p}(\Omega)$. We have that \begin{align*}
    &\langle (-\Delta)_{p,w}^su-(\Delta)_{p,w}^sv,u-v\rangle\\
    &=\left(-\operatorname{div}^s\left(w(x)|\nabla^su|^{p-2}\nabla^2u\right)-\operatorname{div}^s\left(w(x)|\nabla^sv|^{p-2}\nabla^2v\right)\right)(u-v)\\
    &=\int_{\R^n}w(x)|\nabla^su|^{p-2}\nabla^su\cdot (\nabla^su-\nabla^sv)\,dx-\int_{\R^n}w(x)|\nabla^sv|^{p-2}\nabla^sv\cdot (\nabla^su-\nabla^sv)\,dx\\
    &=\int_{\R^n}w(x)\left(|\nabla^su|^{p-2}\nabla^su-|\nabla^sv|^{p-2}\nabla^sv\right)\cdot \left(\nabla^su-\nabla^sv\right)\,dx,
\end{align*}
For vectors $x,y\in \R^n$, we have the identity \cite[Lemma A.0.5]{Ireneo97} $$\left(|x|^{p-2}x-|y|^{p-2}y\right)\cdot (x-y)\geq \begin{cases}
    c_p|x-y|^p,\,p\geq 2,\\
    c_p\frac{|x-y|^2}{\left(|x|+|y|\right)^{2-p}},\ 1<p<2.
\end{cases},$$ for some positive constant $c_p$. Thus \begin{align*}
    \langle (-\Delta)_{p,w}^su-(\Delta)_{p,w}^sv,u-v\rangle\geq \begin{cases}
    2^{2-p}\norm{\nabla^s(u-v)}_{L^p_w(\R^n;\R^n)}^p,\,p\geq 2\\
    (p-1)\int_{\R^n}w(x)\frac{|\nabla^su-\nabla^sv|^2}{\left(|\nabla^su|+|\nabla^sv|\right)^{2-p}}\,dx,\,1<p<2,
    \end{cases}.
\end{align*}
Observe that from H\"older's inequality: \begin{align*}
    &\norm{\nabla^su-\nabla^sv}_{L^p_w(\R^n;\R^n)}^p\\
    &\leq \left(\int_{\R^n}\frac{|\nabla^su-\nabla^sv|^2}{\left(|\nabla^su|+|\nabla^sv|\right)^{2-p}}w(x)\,dx\right)^{p/2}\left(\int_{\R^n}\left(|\nabla^su|+|\nabla^sv|\right)^pw(x)\,dx\right)^{(2-p)/2}
,\end{align*} so 
\begin{align*}&\langle (-\Delta)_{p,w}^su-(\Delta)_{p,w}^sv,u-v\rangle\\
&\geq \begin{cases}
    2^{2-p}\norm{\nabla^s(u-v)}_{L^p_w(\R^n;\R^n)}^p,\,p\geq 2\\
    (p-1)\norm{\nabla^su-\nabla^sv}_{L^p_w(\R^n;\R^n)}^2\left(\norm{\nabla^su}_{L^p_w(\R^n;\R^n)}+\norm{\nabla^sv}_{L^p_w(\R^n;\R^n)}\right)^{p-2}\,1<p<2,
    \end{cases}.\end{align*}
and this establish that the operator $(-\Delta)_{p,w}^s$ is monotone.

Finally, to prove the continuity of the operator, we take $(u_n)\subset X^{s,p}_{0,w}(\Omega)$ converging strongly to some function $u\in X^{s,p}_{0,w}(\Omega)$. Now, let $v\in X^{s,p}_{0,w}(\Omega)$. We have that \begin{align*}
    &\langle (-\Delta)_{p,w}^su_n-(-\Delta)_{p,w}^su,v\rangle=\int_{\R^n}w(x)\left(|\nabla^su_n|^{p-2}\nabla^su_n-|\nabla^su|^{p-2}\nabla^su\right)\cdot \nabla^sv\,dx\\
    &\leq \left(\int_{\R^n}w(x)\left(|\nabla^su_n|^{p-2}\nabla^su_n-|\nabla^su|^{p-2}\nabla^su\right)^{p'}\,dx\right)^{1/p'}\norm{\nabla^sv}_{L^p_w}^.
\end{align*}
Since $u_n\to u$ in $X^{s,p}_{0,w}(\R^n)$, there is a (non relabeled) subsequence such that $\nabla^su_n(x)\to \nabla^su(x)$ a.e. for $x$, and a function $h\in L^p_w(\R^n)$ such that $|\nabla^su(x)|\leq h(x)$ a.e. for $x$. Hence, we can conclude that $$\lim_{n\to\infty}\int_{\R^n}w(x)\left(|\nabla^su_n|^{p-2}\nabla^su_n-|\nabla^su|^{p-2}\nabla^su\right)^{p'}\,dx=0,$$ which shows that $$\lim_{n\to\infty}\langle (-\Delta)_{p,w}^su_n,v\rangle=\langle (-\Delta)_{p,w}^su,v\rangle,\,\forall v\in X_{0,w}^{s,p}(\Omega).$$
Hence, by Theorem \ref{BW}, there exists a weak solution $u\in X_{0,w}^{s,p}(\Omega)$ to the problem \ref{eqn:P}. Moreover, if there exists $u_1,u_2\in X_{0,w}^{s,p}(\Omega)$ two different weak solutions to the problem, the function $u:=u_1-u_2\in X_{0,w}^{s,p}(\Omega)$ solves $$\int_{\R^n}w(x)|\nabla^su|^{p-2}\nabla^su\cdot \nabla^sv\,dx=0,\,\forall v\in X_{0,w}^{s,p}(\Omega).$$ In particular, for the choice $v=u$, we have   $$\norm{\nabla^su}_{L^p_w}^p=\int_{\R^n}w(x)|\nabla^su|^{p}\,dx=0,$$ which implies $u=0$ a.e., and hence the weak solution is unique.
\qed

In the particular case $p=2$, since we are dealing with Hilbert spaces, we can show existence and uniqueness for a more general problem with degenerate ellipticity:
$$\begin{cases}
    -\operatorname{div}^s\left(A(x)\nabla^su\right)&=f,\,\Omega,\\
    \hfill u&=0,\,\Omega^c.
\end{cases},$$ where $A(x):\Omega\to \R^{n\times n}$ is a symmetric matrix satisfying the following degenerate ellipticity condition: there exists two positive constants $c_1,c_2>0$ such that for every $x\in \Omega$ and $\xi\in  \R^n$, $$c_1w(x)|\xi|^2\leq A(x)\xi\cdot \xi\leq c_2 w(x)|\xi|^2,$$ where $w$ is a weight in $A_2$. We can also prescribe the values of $u$ outside $\Omega$ to a function $g\in X^{s,2}_{w}(\R^n)$, and hence we will look to weak solutions on the space $X^{s,p}_{g,w}(\Omega)$ defined as $$X^{s,p}_{g,w}(\Omega):=\{u\in X^{s,p}_w(\R^n): u=g,\,\text{a.e. in}\, \Omega^c\}.$$
\begin{Teor}
    Let $w\in A_2$ and $A(x):\Omega\to \R^{n\times n}$ a symmetric matrix satisfying that there exists two positive constants $c_1,c_2>0$ such that for every $x\in \Omega$ and $\xi\in  \R^n$, $$c_1w(x)|\xi|^2\leq A(x)\xi\cdot \xi\leq c_2 w(x)|\xi|^2.$$ For $f\in H^{-s,2}_{w^{-1}}(\Omega)$, $g\in X^{s,2}_{w}(\R^n)$, there exists an unique $u\in X^{s,2}_{g,w}(\Omega)$ such that $$\int_{\R^n}A(x)\nabla^su\cdot \nabla^sv\,dx=\int_\Omega fv\,dx,\,\forall v\in X^{s,2}_{0,w}(\Omega).$$
\end{Teor}
\noindent{}\textbf{Proof:} The proof follows as a typical application of the Lax-Milgram lemma. Let us suppose that $u$ is a weak solution, and we define $\tilde{u}:=u-g\in X_{0,w}^{s,2}(\Omega)$. By linearity of the fractional operators we have that $\nabla^su=\nabla^s\tilde{u}+\nabla^sg$, and hence the problem in its variational form reads as $$\int_{\R^n}A(x)\nabla^s\tilde{u}\cdot \nabla^sv\,dx=\int_\Omega fv\,dx-\int_{\R^n}A(x)\nabla^sg\cdot \nabla^sv\,dx,\,v\in X^{s,2}_{0,w}(\Omega).$$ We define the bilinear form $a(\tilde{u},v):X^{s,2}_{0,w}(\Omega)\times X^{s,2}_{0,w}(\Omega)\to \R$ as $$a(\tilde{u},v):=\int_{\R^n}A(x)\nabla^s\tilde{u}\cdot \nabla^sv\,dx.$$ Since $A$ is a symmetric matrix, the bilinear form $a$ is symmetric. Let $v=\tilde{u}$, then \begin{align*}
    a(\tilde{u},\tilde{u})&=\int_{\R^n}A(x)\nabla^s\tilde{u}\cdot \nabla^s\tilde{u}\,dx\geq c_1\int_{\R^n}|\nabla^s\tilde{u}|^2w(x)\,dx=c_1\norm{\nabla^s\tilde{u}}_{L^2_w(\R^n;\R^n)}^2\\
    &\geq c'_1\norm{u}_{X^{s,2}_w(\R^n)}^2,
\end{align*} so $a$ is coercive. The continuity follows by the fact that we can rewrite $$A(x)\nabla^s\tilde{u}\cdot \nabla^sv=\left((\nabla^s\tilde u)^t\left(A(x)^{1/2}\right)^t\right)^tA(x)^{1/2}\nabla^sv,$$ and hence by Cauchy-Schwarz we have that \begin{align*}
    |a(\tilde{u},v)|\leq \left(\int_{\R^n}|A(x)||\nabla^s\tilde{u}|^2\,dx\right)^{1/2}\left(\int_{\R^n}|A(x)||\nabla^sv|^2\,dx\right)^{1/2},
\end{align*} where $|A(x)|$ must be seen as the Frobenius norm of the matrix $A(x)$, and by the hypothesis over $A(x)$ we have that $$|A(x)|\leq Cw(x),\,x\in \Omega,$$ for a positive constant $C$ independent of $x$. Thus, we have that \begin{align*}|a(\tilde{u},v)|&\leq C\left(\int_{\R^n}|\nabla^s\tilde{u}|^2w(x)\,dx\right)^{1/2}\left(\int_{\R^n}|\nabla^sv|^2w(x)\,dx\right)^{1/2}\\
&=C\norm{\nabla^s\tilde{u}}_{L_w^2(\R^n;\R^n)}\norm{\nabla^sv}_{L_w^2(\R^n;\R^n)}<\infty,\end{align*} so the bilinear form $a$ is continuous. Note that the functional $v\mapsto \int_\Omega fv$ is continuous in $X^{s,2}_{0,w}(\Omega)$ since \begin{align*}
    \int_\Omega |fv|\,dx&=\int_\Omega|f|w(x)^{-1/2}
    |v|w(x)^{1/2}\,dx\leq \left(\int_\Omega |f|^2w(x)^{-1}\,dx\right)^{1/2}\left(\int_\Omega |v|^2w(x)\,dx\right)^{1/2}\\
    &=\norm{f}_{L^2_{w^{-1}}(\Omega)}\norm{v}_{L^2_{w(\Omega)}}<\infty.
\end{align*} Hence, by an analogous argument as for the continuity of $a$ we can prove that $$v\mapsto \int_\Omega fv\,dx-\int_{\R^n}A(x)\nabla^sg\cdot \nabla^sv\,dx,$$ is continuous, so by Lax-Milgram's Theorem there exists a function $\tilde{u}\in X^{s,2}_{0,w}(\Omega)$ such that $$a(\tilde{u},v)=v\mapsto \int_\Omega fv\,dx-\int_{\R^n}A(x)\nabla^sg\cdot \nabla^sv\,dx,$$ for every $v\in X^{s,2}_{0,w}(\Omega)$. Since $u=\tilde{u}+g\in X^{s,2}_{g,w}(\Omega)$, $u$ is a weak solution for the problem \ref{eqn:P}. \\
For uniqueness, we take $u_1,u_2\in X_{g,w}^{s,2}(\Omega)$ weak solutions to our problem. We define $u:=u_1-u_2\in X_{0,w}^{s,2}(\Omega)$. By linearity of the fractional operators we have that $$-\operatorname{div}^s\left(A(x)\nabla^su\right)=-\operatorname{div}^s\left(A(x)\nabla^su_1\right)+\operatorname{div}^s\left(A(x)\nabla^su_2\right)=0,$$ and hence $$\int_{\R^n}A(x)\nabla^su\cdot \nabla^s v\,dx=0,\,\forall v\in X_{0,w}^{s,p}(\Omega).$$ In particular, for $v=u$ it yields that $$\int_{\R^n}A(x)\nabla^su\cdot \nabla^su\,dx=0\implies 0\leq c_1\int_{\R^n}|\nabla^su|^2w(x)\,dx\leq 0,$$ and together with the fact that $u$ is $0$ outside $\Omega$ we have that $u=0$ a.e. in $\R^n$, i.e., $u_1=u_2$ a.e. in $\R^n$. \qed

For a subsequent work, we would like to study the regularity of the weighted fractional $p$-Laplacian. In the unweighted Hilbertian case, the operator $Lu=\operatorname{div}\left(A(x)\nabla u\right)$ with $A$ satisfying the degenerate ellipticity condition was studied by E. Fabes, C. Kenig, and R. Serapioni in \cite{Fabes1982}, where they proved local H\"older regularity. Further results has been established more recently by A. Claros in \cite{Claros2025}.
\section*{Statements and Declarations}
This work was supported by {\it Agencia Estatal de Investigación} (Spain) through grant PID2023-151823NB-I00 and {\it Junta de Comunidades de Castilla-La Mancha} (Spain) through grant SBPLY/23/180225/000023. G.G.-S. is supported by a Doctoral Fellowship by \textit{Universidad de Castilla-La Mancha} \text{2024-UNIVERS-12844-404}.
\section*{Acknowledgements}
The author would like to thank to his advisor José Carlos Bellido for the useful suggestions on a preliminary
version of the manuscript and P.M. Campos for several fruitful discussions on the subject during his visit to \textit{Universidad de Castilla-La Mancha} in February of 2025.

\appendix


\section*{Conflicts of interest}

The author declare that there are no conflicts of interest regarding the publication of this paper.

\addcontentsline{toc}{section}{References}
\bibliographystyle{plain}

\end{document}